\documentclass[11pt]{article}
\usepackage{latexsym}
\usepackage{amssymb,amsbsy,amsmath,amsfonts,amssymb,amscd}
\usepackage{epsfig, graphicx}
\usepackage{colortbl}
\newcommand{\fer}[1]{(\ref{#1})}
\newcommand {\wt}[1] {{\widetilde #1}}

\newcommand{\commentout}[1]{}
\newcommand{\R}{\mathbb{R}}

\newcommand{\ep}{\varepsilon}
\newcommand {\al} {\alpha}
\newcommand {\e}  {\varepsilon}

\newcommand {\lb} {\lambda}
\newcommand {\Chi} {{\bf \raise 2pt \hbox{$\chi$}} }

\newcommand {\cac} { {\mathcal C} }

\newcommand {\f}   {\frac}
\newcommand {\p}   {\partial}
\newcommand{\beq}{\begin{equation}}
\newcommand{\beqa} {\begin{array}{rl}}
\newcommand{\eeq}{\end{equation}}
\newcommand{\eeqa}{\end{array}}
\newtheorem{theorem}{Theorem}[section]
\newtheorem{lemma}[theorem]{Lemma}

\newtheorem{proposition}[theorem]{Proposition}

\newcommand{\qed}{{ \hfill
                       {\unskip\kern 6pt\penalty 500
                       \raise -2pt\hbox{\vrule\vbox to 6pt{\hrule width 6pt
                       \vfill\hrule}\vrule} \par}   }}

\newcommand{\FFT}{{\mathsf{F}}}
\newcommand{\IFFT}{{\mathsf{F}^{\ast}}}
\title{\Large \bf Numerical Solution of an Inverse Problem in Size-Structured Population Dynamics}

\author{Marie Doumic \thanks{{ D\'epartement de Math\'ematiques et Applications,
\'Ecole Normale Sup\'erieure, INRIA projet BANG, 
            45 rue d'Ulm, F~75230 Paris cedex 05, France;
email: doumic@dma.ens.fr}} \footnotemark[3]
\and Beno\^ \i t Perthame\thanks{{ Universit\'e Pierre et Marie Curie-Paris 6, UMR 7598 LJLL, BC187, 4, place Jussieu,  F-75252 Paris cedex 5, and Institut Universitaire de France; email: perthame@ann.jussieu.fr}} \thanks{INRIA Rocquencourt, projet BANG, Domaine de Voluceau, BP 105, 781153 Rocquencourt, France; emails: marie.doumic@inria.fr, benoit.perthame@inria.fr} 
\and Jorge Zubelli\thanks{IMPA, Est. D. Castorina 110, Rio de Janeiro, RJ 22460-320 Brazil;
email: zubelli@impa.br}}

\date{\today}

\begin{document}
\maketitle
\pagestyle{plain}
\pagenumbering{arabic}

\begin{abstract}
We consider a size-structured model for cell division and address the question of determining the division (birth) rate from the measured stable size distribution of the population. We propose a new regularization technique based on a filtering approach. We prove convergence of the algorithm 
and  validate the theoretical results by implementing numerical simulations, based on classical techniques. We compare the results  for direct and inverse problems, for the filtering method and for the quasi-reversibility method proposed in \cite{PZ}. \\

\end{abstract}

\section{Introduction}

The use of size-structured models to describe biological systems has attracted the
interest of many authors and has a long standing tradition.  In particular, the use of
 size structures was very well documented and compared to experiments
 in the 70's. This led to the survey book \cite{MD} and  subsequent mathematical analysis (see also the references in \cite{BP}).
Needless to say, in such models it is crucial for the analysis, computer simulation and prediction to calibrate the
corresponding model parameters so as to obtain good quantitative results.
Indeed, in the inverse problem literature, a number of authors
have addressed the calibration of certain structured population models. See for example~\cite{EnRunSche,GOP,R89,PR91} and references therein.

In this article, we consider theoretical and numerical aspects of the inverse problem of determining the division rate coefficient $B=B(x)$ in the following specific size-structured model for cell division:
\beq
\label{eq:gencell}
\left \{ \begin{array}{l}
 \f{\p}{\p t} n(t,x) +	\f{\p}{\p x} n(t,x)  + B(x) n(t,x)  = 4  B(2 x) n(t,2 x), \qquad x \geqslant0,\, t \geqslant0,
\\
\\
n(t,x=0)=0, t> 0,
\\
\\
n(0,x) = n^{0}(x) \ge 0 .
\end{array} \right.
\eeq
Here, the cell density is represented by $n(t,x)$ at time $t$ and  size $x$. The division rate $B$ expresses the division of cells of size $2x$ into two cells of size $x.$

By making use of flux cytometry technologies for instance, it is possible to determine cell populations with certain properties as protein content on a large scale of tenths of thousands of cells. In other applications, like coagulation fragmentation equation \cite{CarrilloGoudon, ColletGoudon, ColletGoudonVasseur, Laurencot, LaurencotMischler}, or prion aggregation and fragmentation \cite{CP,MJN99,GPW2006}, similar equations arise, and much less is known on aggregate size repartition.  The division rate $B(x)$, on the contrary, is not directly measurable.

The long time behavior of solutions is well known. Indeed, it was proved in \cite{PR, M} that under fairly general conditions on the coefficients, there is a unique solution $(N,\lb_0)$  to the following eigenvalue problem
\beq
\label{eq:celldiv1}
\left \{ \begin{array}{l}
 \f{\p}{\p x} N + (\lb_0+B(x)) N =4  B(2 x) N(2 x), \qquad x \geqslant0,
\\
\\
N(x=0)=0,
\\
\\
N(x)>0 \; \text{ for } x>0, \qquad \int_0^\infty N(x)dx =1, 
\end{array} \right.
\eeq
where $\lb_0>0$ and $N e^{\mu x} \in L^\infty \cap L^1$ for all $\mu < \lb_0$.

It was shown in \cite{MMP,PR} 
$$
n(t,x) e^{-\lb_0 t} {\;}_{\overrightarrow{\; t \rightarrow \infty \;}}\;  m_0 N(x), \quad \mbox{ in }\;L^1(\R_+,\phi(x)dx),
$$
where the weight $\phi$ is the unique solution to the adjoint problem 
\beq
\label{eq:celldiv1:adj}
\left \{ \begin{array}{l}
 -\f{\p}{\p x} \phi + (\lb_0+B(x)) \phi =  2B(x)\phi(\f{x}{2}), \qquad x \geqslant0,
\\
\\
\phi(x)>0, \qquad \int_0^\infty \phi(x)N(x)dx =1.
\end{array} \right.
\eeq
In other words, $\lb_0$ is the growth rate of such a system and is usually called ``Malthus parameter'' in population biology. From \cite{MMP,PR, BP} we also know  that $\lb_0$ is related to $N$ by the relation
\beq\label{eq:lbn}
\lb_{0}= \f{\int_0^{\infty} N dx}{\int_0^{\infty} x N dx}.
\eeq

The question we address here is the following: How can we estimate the division rate $B$ from the knowledge of the steady dynamics $N$ and  $\lb_0$ ? 
The inverse problem thus consists of finding $B$ a solution to 
\beq\label{eq:exact:inverse}
 4B(2x)N(2x) - B(x) N(x) =L(x):= \f{\p}{\p x} N(x) + \lb_0 N(x) ,\quad x\geqslant 0,
\eeq
assuming that $(N,\lb_0)$ is known, or thanks to \fer{eq:lbn} that $N$ is known. 
As seen in \cite{PZ}, this problem is well-posed if $N$ satisfies strong regularity properties such as $\f{\p}{\p x} N({x})\in L^p (\R_+)$ for some $p\geqslant 1.$ 

However, in practical applications we have only an \emph{approximate} knowledge of $(N,\lb_0),$ given by noisy data $(N_\ep,\lb_\ep),$ with $N_\ep\in L^2_+(\R_+)$ for instance.~\footnote{Actually, our knowledge of $\lambda_0$ is presumably an order of precision higher than that of $N$, since the rate $\lambda_0$ can be estimated independently by means of time information.} 
This means that we have no way of controlling $\f{\p}{\p x} N_\ep,$ so we cannot control the precision of a solution $B_\ep$ to problem (\ref{eq:exact:inverse}) when a perturbed $N_\ep$ replaces  $N$. Furthermore, it is not even clear whether such a $B_\e$ exists.

\

The question we focus on is then: How to approximate the problem (\ref{eq:exact:inverse}) in order to get a solution $B_\ep$ as close as possible to the exact division rate $B$?

\

We remark that, in the context of noisy data, the inverse problem under consideration is ill-posed~\cite{PZ} and thus regularization
would be required. 
A natural tool to be invoked from the inverse problem literature would be some kind of Tikhonov regularization method \cite{BaLe,EnHaNe}. However, this would lead to computationally intensive problems.
Indeed, for each forward problem evaluation a dilation-differential equation of the form (\ref{eq:celldiv1})
would have to be solved.

In \cite{PZ}, two of the present authors proposed a method of regularization consisting in the solution of the following approximate problem:
\begin{equation*}
\left\{
\begin{array}{l}
\al \f{\p}{\p y} (B_{\e,\al} N_\e) + 4B_{\e,\al}(y) N_\e(y) = B_{\e,\al}\big(\f y 2\big) N_\e\big(\f y 2\big)+ \lb_0 N_\ep\big(\f y 2\big)+2 \f{\p}{\p y}\biggl( N_\e\big(\f y 2\big)\biggr), \qquad y>0,
\\
\\
(B_{\e,\al} N_\e)(0)=0,
\end{array} \right.
\end{equation*}
where $\alpha$ is a regularizing parameter. It was shown that a convergence rate of order $\sqrt{\ep}$ could be obtained, for $\alpha=O(\sqrt{\ep}),$ where $\ep$ is the error on the data $N$ in an appropriate norm.

The above method of the solution to the inverse problem will be called {\em quasi reversibility} in accordance with the general spirit of the terminology of \cite{Lattes,LattesLions}. 
The main goal of this work is to investigate the numerics of such approach, to consider an  alternative technique based on filtering ideas and to compare the performance of the
different methods. The alternative technique is also analyzed from the theoretical point of view
and estimates are presented.

In this work, we have modified slightly the original regularization equation by writing $\lb_{\ep,\alpha}$ instead of $\lb_0$ for the reasons we shall explain in the sequel.
Thus, we work with 
\beq
\left\{
\begin{array}{l}
\al \f{\p}{\p y} (B_{\e,\al} N_\e) + 4B_{\e,\al}(y) N_\e(y) = B_{\e,\al}\big(\f y 2\big) N_\e\big(\f y 2\big)+ \lb_{\ep,\alpha} N_\ep\big(\f y 2\big)+2 \f{\p}{\p y}\biggl( N_\e\big(\f y 2\big)\biggr), \qquad y>0,
\\
\\
(B_{\e,\al} N_\e)(0)=0.
\end{array} \right.
\label{eq:invfull}
\eeq
Indeed, in order to conserve regularity properties of the solution $H=BN$ to the inverse problem, we want it to be both in $L^1 (\R_+)$ and in $L^1 (\R_+,xdx)$ 
in order to express that both the total number of cells and the total biomass are finite. 
Hence, formal integration of Equation (\ref{eq:invfull}) gives
\beq\label{eq:conserv1:PZ}
\lb_{\ep,\alpha} \int_0^{\infty} N_\ep dx=\int_0^{\infty} B_{\ep,\alpha} N_\ep dx,
\eeq
and integration against the weight $x$ gives
\beq\label{eq:integ:x}
-\alpha \int_0^{\infty} B_{\ep,\alpha} N_\ep dx = 4 \lb_{\ep,\alpha} \int_0^{\infty} x N_\ep dx - 4\int_0^{\infty} N_\ep dx.
\eeq 
Hence, we have to choose, according to the eigenvalue theory:
\beq\label{eq:conserv2:PZ}
\lb_{\ep,\alpha}= \f{\int_0^{\infty} N_\ep dx}{\int_0^{\infty} x N_\ep dx + \f{\alpha}{4}\int_0^{\infty} N_\ep dx}.
\eeq
The choice of $\lb_{\ep,\alpha}$ can be understood as a compatibility condition when $\alpha > 0$ and for
$\alpha = 0$ it tells us that $(N,\lb_0)$ is overdetermined data for the inverse problem. 
Therefore, if we have {\em a priori} knowledge on $\lambda_0$, we could verify its distance to $\lambda_{\epsilon,\alpha}$ as a way of checking the error of the inverse problem solution.

\

The plan of this work is the following: 
In Section~\ref{sec:filter}, we propose yet another method to regularize the inverse problem, and obtain a convergence rate. The convergence rate turns out to be as good as the one in \cite{PZ}. In Section~\ref{sec:scheme} we give a numerical method to solve it, and in Section~\ref{sec:result} we show some numerical simulations so as to compare the accuracy of the different methods.

\section{Regularization by Filtering}
\label{sec:filter}
\subsection{Filtering approach}

Taking a closer look at Equation~(\ref{eq:exact:inverse}), we see that all the difficulties come from the differential term  $\f{\p}{\p x} N.$ In \cite{PZ}, the choice was to add an equivalent derivative $\alpha \f{\p}{\p x} (B\,N)$ to the equation; here on the contrary, we choose to regularize it by a convolution method.

For $\alpha>0$,  we use the notation 
\beq\label{eq:defrhoalpha}
\rho_\alpha(x) = \f{1}{\alpha}\rho(\f{x}{\alpha}),\qquad \rho\in \cac_c^\infty(\R),\quad\int_0^{\infty}\rho(x)\,dx=1,\quad\rho\geqslant 0,\quad \mathrm{Supp}(\rho)\subset [0,1], 
\eeq
and we replace in (\ref{eq:exact:inverse}) the term $\f{\p}{\p x} N_\ep + \lb_0 N_\ep$ by
$$\bigl(\f{\p}{\p x} N_\ep + \lb_{\ep,\alpha} N_\ep\bigr)* \rho_\alpha (x)=
N_\ep * \bigl(\f{\p}{\p x} \rho_\alpha + \lb_{\ep,\alpha} \rho_\alpha \bigr)(x)=\int_0^{\infty} N_\ep (x') \bigl(\f{\p}{\p x} \rho_\alpha +\lb_{\ep,\alpha} \rho_\alpha\bigr)(x-x') dx'.$$ 
We now use the notation
$$N_{\ep,\alpha}=N_\ep * \rho_\alpha.$$
In this way, we obtain a smooth term in $L^2(\R_+)$. Furthermore, $N_{\ep,\alpha}$   converges  
to $N_\ep$ in $L^2(\R_+)$ when $\alpha$ tends to zero. 
We now have to consider the following problem: 
\\
Find $B_{\ep,\alpha}$ solution of
\beq
\label{eq:inverse:filter}
4B_{\ep,\alpha}(2x)N_{\ep,\alpha}(2x)  + B_{\ep,\alpha}(x) N_{\ep,\alpha}(x)= \f{\p}{\p x} N_{\ep,\alpha} + \lb_{\ep,\alpha}  N_{\ep,\alpha}(x),\quad x\geqslant 0 \mbox{ .}
\eeq
As in Equation~(\ref{eq:invfull}), for the quasi-reversibility method, we need to choose $\lb_{\ep,\alpha}$ appropriately.
Indeed, we perform the same manipulations leading to Equation~(\ref{eq:conserv2:PZ}) to get
\beq\label{eq:lambdaep}
\lb_{\ep,\alpha} = \f{\int_0^{\infty} N_{\ep,\alpha}(x) dx}{\int_0^{\infty} x N_{\ep,\alpha} (x) dx}.
\eeq

By Theorem~\ref{th:exist} (see the Appendix), we know that the problem in Equation~(\ref{eq:inverse:filter}) has a unique solution $B_{\e,\alpha} \in L^2(\R_+,N_{\ep,\alpha}^2 dx).$ 

\subsection{Estimates for the filtering approach}
The main result of this section establishes an estimate for the regularization of the 
inverse problem by means of the filtering method described above. 
\begin{theorem}\label{th:estim:filter} Suppose that $N\in H^2 (\R_+)$ and $B\in L^\infty (\R_+),$ $B\geqslant 0$ verify  (\ref{eq:celldiv1}).
Let $\ep>0$ and $N_\ep \in L^2(\R_+),$ $N_\ep (x) >0$ for $x>0,$ such that $$||N_\ep - N||_{L^2(\R_+)}\leqslant \ep ||N||_{L^2 (\R+)}.$$
Let $B_{\ep,\alpha} \in L^2(\R_+,N_{\ep,\alpha}^2dx)$ be the unique solution of (\ref{eq:defrhoalpha}) and (\ref{eq:inverse:filter}).
We have the following estimate:

\beq\label{ineq:estim:filter}
||B_{\ep,\alpha}-B||_{L^2(N_{\ep}^2 dx)} \leqslant C (\alpha+|\lb_{\ep,\alpha} - \lb_0|) ||N||_{H^2(\R_+)} + \f{C}{\alpha} ||N_{\ep,\alpha} - N||_{L^2(\R_+)},
\eeq
where $C$ is a constant depending only on $||B||_{L^\infty},$ $||B_{\ep,\alpha}||_{L^\infty}$ and the regularizing function $\rho.$
\end{theorem}
This theorem relies on a first estimate.

\begin{proposition}\label{prop:estim:filter}
Using the same notations as in Theorem~\ref{th:estim:filter}, we have
\beq\label{ineq:estim:filter2}
||B_{\ep,\alpha}N_{\ep,\alpha}-B\,N||^2_{L^2( dx)} \leqslant C \left( 1+\lb_0^2 \right)\left(1+ \f{1}{\alpha^2}\right)  ||N_{\ep} - N||^2_{L^2(dx)} +C(\alpha^2 + |\lb_{\ep,\alpha} - \lb_0|^2) ||N||^2_{H^2(\R_+)} \mbox{ ,} 
\eeq 
where $C$ depends only  on the regularizing function $\rho.$
\end{proposition}
{\bf{Proof of Prop.~\ref{prop:estim:filter}:}}
Denote by $Q=B_{\ep,\alpha} N_{\ep,\alpha}-B N,$ $R=N_{\ep,\alpha} - N$ and $\delta=\lb_{\ep,\alpha} - \lb_0.$ 
 From Equations~(\ref{eq:celldiv1}) and (\ref{eq:inverse:filter}), $Q$ verifies:
\beq\label{eq:Q}
\left \{ \begin{array}{l}
 \f{\p}{\p x} R (x) +\lb_0 R(x) + \delta N_{\ep,\alpha} (x) + Q(x) =4 Q(2 x), \qquad x \geqslant 0,
\\
\\
Q(x=0)=0.
\end{array} \right.
\eeq
(Since $N_{\ep,\alpha} \in H^1 (\R_+),$ the definition of $Q(x=0)$ is not ambiguous.)
Multiplying this equation by $Q(2x)$ and integrating on the interval $(0,y)$ yields
$$
4 \int\limits_0^y Q (2x)^2 dx =
\int\limits_0^y Q(2x) \f{\p}{\p x} R (x) dx
+ \lb_0 \int\limits_0^y Q(2x) R(x) dx
$$
$$+ \,\delta \int\limits_0^y Q(2x) N_{\ep,\alpha} (x) dx
+\int\limits_0^y Q(2x) Q(x)dx.
$$
 From the Cauchy-Schwarz inequality, after the change of variables $x\rightarrow 2x,$ we have 
$$4\int\limits_0^y Q(2x)^2 dx \leqslant \f{1}{2}\int\limits_0^y \biggl( \f{\p}{\p x} R\biggr)^2 (x) dx
+\f{1}{2}\int\limits_0^y Q (2x)^2 dx
+\f{\lb_0}{2}\int\limits_0^y C R(x)^2 dx + \f{\lb_0}{2} \int\limits_0^y \f{Q(2x)^2}{C} dx 
$$
$$+\,\f{|\delta|^2}{2} \int\limits_0^y N_{\ep,\alpha} (x)^2 dx + \f{1}{2}\int\limits_0^y Q (2x)^2dx
+\f{1}{2}\int\limits_0^y Q(2x)^2 dx
+\int\limits_0^{\f{y}{2}} Q(2x)^2dx.
$$
We take, for instance,
 $C={\lb_0}$. We obtain
\beq ||B_{\ep,\alpha} N_{\ep,\alpha}-B N||^2_{L^2} \leqslant 
 ||N_\ep*\f{\p}{\p x} \rho_\alpha - \f{\p}{\p x} N||^2_{L^2}
+ \lb_0^2 ||N_{\ep}*\rho_\alpha-N||^2_{L^2}
+{|\lb_{\ep,\alpha} -\lb_0|^2} ||N_\ep*\rho_\alpha||^2_{L^2}.
\label{ineq:prop}
\eeq 
The last two terms of this inequality are easy to estimate, writing
$$ ||N_{\ep}*\rho_\alpha-N||_{L^2} \leqslant ||N_{\ep}*\rho_\alpha-N*\rho_\alpha||_{L^2}  + ||N*\rho_\alpha-N||_{L^2} \leqslant C \left(||N_\ep-N||_{L^2} + \alpha ||N||_{H^1}\right), $$
and
$$ ||N_\ep*\rho_\alpha||_{L^2}\leqslant C ||N||_{L^2}.$$
It remains to evaluate the first term on the right-hand side of inequality (\ref{ineq:prop}). 
We write 
$$||N_\ep*\f{\p}{\p x} \rho_\alpha - \f{\p}{\p x} N||^2_{L^2} \leqslant 2 ||N_\ep*\f{\p}{\p x} \rho_\alpha - N*\f{\p}{\p x} \rho_\alpha||^2_{L^2}
+
2 ||N *\f{\p}{\p x} \rho_\alpha - \f{\p}{\p x} N||^2_{L^2}.$$
By a convolution estimate we evaluate the first term as
$$||N_\ep*\f{\p}{\p x} \rho_\alpha - N*\f{\p}{\p x} \rho_\alpha||^2_{L^2(\R_+,dx)} \leqslant
||N_\ep - N||^2_{L^2(\R_+,dx)} ||\f{\p}{\p x} \rho_\alpha ||^2_{L^1} \mbox{ .}$$
Since $\int_0^{\infty} | \f{\p}{\p x} \rho_\alpha (x) | dx= \f{1}{\alpha}\int_0^{\infty} | \f{\p}{\p x} \rho (y) | dy ,$ we have
$$||N_\ep*\f{\p}{\p x} \rho_\alpha - \f{\p}{\p x} N||^2_{L^2(\R_+,dx)} \leqslant 
\f{C(\rho)}{\alpha^2}||N_\ep - N||^2_{L^2(\R_+,dx)}
+
2 ||N *\f{\p}{\p x} \rho_\alpha - \f{\p}{\p x} N||^2_{L^2(\R_+,dx)} \mbox{ .}$$
To evaluate the last term $||N *\f{\p}{\p x} \rho_\alpha - \f{\p}{\p x} N||^2_{L^2},$
we extend to $\R$ the functions $N$ and $\f{\p}{\p x} N$ by zero and consider their Fourier transforms. We denote $\hat{f}(\xi)$ the Fourier transform of $f\in L^2(\R_+)$ at $\xi$, where $f$ is extended as zero on $\R_-.$ 
We obtain by Fourier analysis
$$||N *\f{\p}{\p x} \rho_\alpha - \f{\p}{\p x} N||^2_{L^2(\R_+,dx)}
=
||i\xi \hat{N} \hat{\rho}_\alpha - i\xi\hat{N}||^2_{L^2(\R_+,dx)}\leqslant \int_{-\infty}^\infty |\hat{N}(\xi)|^2 |\xi|^4 \f{|\hat{\rho}_\alpha (\xi)-1|^2}{|\xi|^2} d\xi. 
$$ 
Using that 
\beq\label{ineq:rhoalpha}
|\f{|\hat{\rho}_\alpha(\xi)-1|^2}{\xi^2}|       
 \leqslant C(\rho) \alpha^2,
\eeq
where $C(\rho)$ only depends on the regularization function $\rho$, we have that 
$$
||N *\f{\p}{\p x} \rho_\alpha - \f{\p}{\p x} N||^2_{L^2(\R_+,dx)} \leqslant C(\rho) \alpha^2 ||N||^2_{H^2(\R_+)}  \mbox{ .}$$
Going back to (\ref{ineq:prop}), this concludes the proof of Proposition \ref{prop:estim:filter}.
\qed 

\

We can now deduce the proof of Theorem \ref{th:estim:filter}. We write:
$$||B_{\ep,\alpha} - B||_{L^2(N_\ep^2 dx)} \leqslant  ||B_{\ep,\alpha} N_\ep - B_{\ep,\alpha}\,N_{\ep,\alpha}||_{L^2(\R_+)} + ||B_{\ep,\alpha}\,N_{\ep,\alpha} - B\, N||_{L^2 (\R_+)}+ ||B\,N - B\, N_\ep||_{L^2 (\R_+)}.
$$ 
Using Proposition~\ref{prop:estim:filter}, and the fact that
$$||N_\ep - N_{\ep,\alpha} ||_{L^2} \leqslant 2 ||N_\ep - N||_{L^2} + \alpha ||N||_{H^1} ,$$
this inequality gives the result. 
\qed

\section{Numerical Solution of the Inverse Problem}
\label{sec:scheme}

This section is concerned with the numerical aspects of the solution of the inverse problem.
In order to do that we start with a description of the solution to the direct one in Subsection~\ref{subsec:direct}.

\subsection{Direct Problem}\label{subsec:direct}

In the direct problem, we assume we know the proliferation rate $B$, we look for $N$ and $\lb_0>0$ solutions of (\ref{eq:celldiv1}).
For this purpose, we solve the time-dependent problem (\ref{eq:gencell}) and look for a steady dynamics. As already said,
this problem is well-posed (see for instance \cite{BP}) and it was proved in \cite{MMP} that solutions grow at an exponential rate towards $\rho N(x)e^{\lb_0 t}$ with  $\rho=\int_0^{\infty} n(0,x) \phi(x)  dx,$ recalling the notation in (\ref{eq:celldiv1:adj}). Furthermore, under more restrictive conditions it was shown in \cite{PR} that there exists constants $\mu>0$ and $C(n^0)>0,$ such that 
$$
||n(t,x) e^{-\lb_0 t} -\rho N(x)||_{L^1(\R_+, \phi (x)dx)} \leqslant C e^{-\mu t}.
$$ 
To solve it numerically, we discretize the problem (\ref{eq:gencell}) along a regular grid, denote by $\Delta t$ the time step and by $\Delta x={L}/{I}$ the spatial step, 
 where $I$ denotes the number of points and $L$ the computational domain length: $x_i=i\Delta x,$ $0\leqslant i\leqslant I.$

We use an upwind finite volume method (cf. \cite{Bouchut, LeVeque, Godlewski})
$$
n_i^k= \f{1}{\Delta x} \int\limits_{x_{i-\f{1}{2}}}^{x_{i+\f{1}{2}}} n(k\Delta t,y)dy,\qquad  \f{1}{\Delta t}\int\limits_0^{\Delta t} n(k\Delta t +s, x_{i+\f{1}{2}}) ds  \approx n_i^k.
$$
For the time discretization, 
we use a marching technique. 
We choose the time step $\Delta t$ so as to satisfy the largest possible CFL stability criteria $\theta:= \f{\Delta t}{\Delta x}=1$.

The numerical scheme is given, for $i=1,...,I,$ by $n_0^k=0$ and 
\beq\label{eq:scheme:direct}
\f{n_i^{k+1} - n_i^k}{\Delta t} + \f{n_i^k - n_{i-1}^k}{\Delta x} 
+
 B_i n_i^{k+1}= B_{2i-1} n_{2i-1}^k + 2 B_{2i} n_{2i}^k+B_{2i+1}n_{2i+1}^k \mbox{ ,}
\eeq
with the convention that $n_j=0$ for $j>I$.  
For stability reasons, we have used an implicit method for the division term in the left hand side and explicit for the right hand side of the equation. The specific form for the right hand side is simply motivated by the need of also dividing cells of odd labels.

According to the power algorithm, we do not keep $n^{k+1}$ from (\ref{eq:scheme:direct}) but 
rather renormalize it as 
$$\wt n^{k+1}= \f{n^k}{\Delta x \sum\limits_{j=1}^{I} n_j^k}.
$$
It is standard, for these positive matrices arising in (\ref{eq:scheme:direct}), that
$$\wt n ^{k+1} \underset{k \to \infty}{\longrightarrow}  N,\qquad \sum\limits_{i=1}^I N_i=1, \qquad N_i >0,
$$
where $N$ is the dominant eigenvector for the problem
$$ \f{N_i -N_{i-1}}{\Delta x} 
+
(\lb_0 + B_i)N_i= B_{2i-1}N_{2i-1} + 2 B_{2i} N_{2i}+B_{2i+1}N_{2i+1}.
$$
One can also find the dominant eigenvalue as 
$$\lb_0=\lim\limits_{k\to\infty} \f{1}{\Delta t} \log\biggl(\f{\sum\limits_{i=1}^I n_i^{k+1}}{\sum\limits_{i=1}^I n_i^k}\biggr).$$
For matrices with one dominant eigenvalue and a corresponding one-dimensional eigenspace, 
it is known that the power algorithm is fast and in fact converges with exponential rate  \cite{Serre}. In practice we can stop the iterations when the relative error  on the normalized quantity
$$\f{1}{\Delta t} \biggl(\sum\limits_{i=1}^I \wt n_i^{k+1} - \sum\limits_{i=1}^I \wt n_i^k\biggr)$$ is small enough, say of the order of $10^{-10}$.

\subsection{Inverse Problem: General Strategy}
\label{sec:strategy}

In the sequel, we denote by $H$ the product $B.N$ and its approximations. 
Indeed, from Equations~(\ref{eq:invfull}) or (\ref{eq:inverse:filter}), we have to search for the product $H=B_{\ep,\alpha}N_\ep$ or $H=B_{\ep,\alpha}N_{\ep,\alpha}$ before computing $B_{\ep,\alpha}.$ In particular, we cannot avoid a loss of information where $N_\ep$ is small, i.e., for $x \approx 0$ or $x\gg 1$. 

The inverse problem  (\ref{eq:exact:inverse}), as well as (\ref{eq:inverse:filter}), can be written as
\beq\label{eq:BN} 
4 H (2x) - H(x) =L(x),
\eeq
with different expressions for $H$ and $L$.  We may think of two possible numerical approaches. 
\\[4mm]
{\bf Strategy 1}. Compute $H(2x)$ from $H(x)$: This means that we re-write Equation  (\ref{eq:BN}) with the new variable  $y=2x$, and arrive at  
\beq\label{eq:BNy} 
4 H (y) - H(\f{y}{2}) =L(\f{y}{2}).
\eeq
The scheme departs from zero, and one deduces the values of $H_i$ step by step, from the knowledge of $H_{j}$ for $j\leqslant i-1$.
\\[2mm]
{\bf Strategy 2}. Compute $H(x)$ from $H(2x)$: The scheme departs from the largest point $x=L$ of our simulation domain. We suppose that for $x\geqslant L$ we have $H(x)=H(L)=0$ (it is relevant since we suppose that $N$ vanishes for $x$ large: see below), and then deduce the smaller values $H_i$ step by step, from the knowledge of $H_{j}$ for $j \geqslant i+1$. 

\

 The two approaches do not necessarily lead to the same result because the continuous equation
\beq\label{eq:H:homogeneous}
4 H(2x) - H(x)=0
\eeq
has infinitely many solutions. This issue is interesting on its own and is related to the construction of wavelets, see \cite{Strang}. It is discussed in Proposition~\ref{prop:pbmH} of the Appendix.

By imposing $H\in L^2 (\R_+),$ we select a unique solution, as shown in Theorem~\ref{th:exist}. The question is then: Which numerical strategy should we use to select the \emph{correct} solution, \emph{i.e.} the one in $L^2(\R_+)$ ?

Among the solutions of Equation (\ref{eq:BN}),  we single out two, defined by the power series:
$$ 
H^{(1)}(x)=\sum\limits_{n=1}^{+\infty} 2^{-2n} L(2^{-n}x)\quad  \mbox{ and } 
H^{(2)} (x)=-\sum\limits_{n=0}^{+\infty} 2^{2n} L(2^n x)
\quad, \quad \forall \;x>0 .
$$
Proposition~\ref{prop:pbmH} shows that for $L\in L^2(\R_+,x^pdx)$, there is a unique solution in $L^2(\R_+,x^p dx),$ given by $H^{(1)}$ if $p<3$ and by $H^{(2)}$ if $p>3$ (and the power series converge in the corresponding spaces).

For $B>0$ smooth and bounded from above and from below, we know that $N$ is smooth and vanishes at $x\approx 0$ and $x\approx \infty,$ and $BN$ inherits these properties. For instance, we know that $H\in L^2 (dx) \cap L^2(x^4 dx).$ 
By uniqueness of a solution in each space, Proposition~\ref{prop:pbmH} implies that $H^{(2)}=H^{(1)},$ or equivalently:
$$
\sum\limits_{n=-\infty}^{+\infty} 2^{2n} L(2^n x)=0 , \quad  \forall\;x\geqslant 0 .
$$
This very particular property cannot be verified at the discrete level.
Hence, the two strategies generally give two different approximations of the same solution of (\ref{eq:BN}).
The first strategy selects an approximation of the solution $H^{(1)}$ whereas the second selects an approximation of the solution $H^{(2)}.$ In the case of a very regular data $N$, then $H^{(2)}$ will perform better around infinity, whereas $H^{(1)}$ will be better around zero. However, if $N$ is a solution of Equation~(\ref{eq:celldiv1}), when we increase the number of points, the two approaches converge to the same solution since $H^{(2)}=H^{(1)}.$ 

Since our simulation domain $[0,L]$ is bounded and contains zero, we prefer the first strategy. This choice is confirmed by all the numerical tests we have performed: the second approach has always lead to a solution exploding around 
zero. However, for the sake of completeness, we also describe the scheme we used for the second approach.


\subsection{Inverse Problem: Filtering Approach}
\label{subsec:filter}

According to strategies 1 and 2, we now present two approaches to handle the numerical solution of the inverse problem regularized with the {\em filtering approach}. Both need to first compute the convolution terms arising in  \fer{eq:inverse:filter}. To do so we first take the Fast Fourier Transform $\FFT$ of $N_\ep$,  multiply it by $i\xi \hat{\rho}_\alpha(\xi),$ and then take the inverse Fast Fourier Transform $\IFFT$.
We choose and define the regularization function $\rho_\alpha$ by its Fourier transform:
$$
\hat{\rho}_\alpha (\xi) =\f{1}{\sqrt{1+\alpha^2\xi^2}}.
$$
This leads us to the numerical approximation
\beq \label{eq:regfft}
\f{\p}{\p x} N_{\ep,\alpha}  \approx dN_\alpha =\mathrm{\IFFT}\biggl(i\xi\hat{\rho}_\alpha (\xi) \mathrm{\FFT}(N_\ep)(\xi)\biggr).
\eeq
We also impose $dN_{\alpha,0}=0$ for compatibility with the continuous equation and further use. 

As mentioned earlier, there are two alternatives, either starting from zero or  coming from infinity. 

\paragraph{The Filtering Approach Starting from Zero (strategy 1).}

We solve Equation~(\ref{eq:inverse:filter}) considered as an equation in the variable $y=2x$, that is to say \fer{eq:BNy}, in order to compute its solution $H^{(1)}(x)$. At the discrete level, we use the notations 
$$
H^f_i\approx B_i N_i, \qquad L^f_i=dN_{\alpha,i} + \lb_{\ep,\alpha} N_i, \quad  L^f_0=0.
$$ 
The discrete version of  \fer{eq:BNy} reads
\beq\label{scheme:filter}
4H^f_{i} = H^f_{\f{i}{2}} + L^f_{\f{i}{2}},\qquad\forall \;0\leqslant i\leqslant I,
\eeq 
and we need to  define the quantities $G_{\f{i}{2}}$. We choose
\beq\label{def:G1/2}
G_{\f{i}{2}}=\left\{\begin{array}{ll}
 G_{\f{i}{2}} \quad & \text{when }i\text{ is even,}
\\ [4mm]
\f{1}{2}\bigl(G_{\f{i-1}{2}} + G_{\f{i+1}{2}}\bigr)\quad & \text{when }i\text{ is odd.}
\end{array}\right.
\eeq 
In particular, we have $H^f_0=0.$

Summing up all the terms in \fer{scheme:filter} for $1\leqslant i \leqslant I,$ we find (with $I$ even to simplify):
$$
4\sum\limits_{i=0}^I H_i^f =2 \sum\limits_{i=0}^{\f{I}{2}} (H^f_i+L^f_i) - \f{1}{2} (H^f_{I/2}+L^f_{I/2}).
$$
Since we have assumed that $N$ has exponential decay for $x\gg 1,$ it follows that
\beq\label{eq:balance1}
\sum\limits_{i=0}^I H^f_i=\sum\limits_{i=1}^{I} L^f_i + E_I,\qquad {\text{with}} \qquad |E_I|\leqslant 2 \sum\limits_{i=\f{I}{2}}^I |H_i|.\eeq
Multiplying (\ref{scheme:filter}) by $x_i$ and summing up again, we find
\beq \label{eq:balance2}
\sum\limits_{i=1}^{\f{I}{2}} x_i L^f_i=F_I,\qquad {\text{with}}\qquad |F_I| \leqslant \sum\limits_{i=\f{I}{2}}^I x_i |H_i|.
\eeq 
As a consequence, we can choose:
\beq \lb_{\ep,\alpha} = -\f{\sum x_i dN_{\alpha,i} }{\sum x_i N_i},\label{def:lbep:f}
\eeq 
as the discrete version of the relations \fer{eq:lbn} or (\ref{eq:lambdaep}).

\paragraph{ The Filtering Approach Starting from Infinity (strategy 2).}

Another method is to discretize the  formulation  (\ref{eq:BN}) in order to compute its solution $H^{(2)}(x)$.
We define the extension $H^f_i=0$ for $i\geqslant I+1$, and for $2\leq i \leq I$, we define by backward iterations
 \beq\label{scheme:infinity}
 H^{f\infty}_i=2H^{f\infty}_{2i}+H^{f\infty}_{2i+1}+H^{f\infty}_{2i-1}-L^f_{i}.
 \eeq
This however does not apply to the indices $i=0,\; 1$ and we set   $H_0^f=\f{L^f_0}{3}=0$ and $H_1^f=4H_2^f -L_1^f$. By summing up all the terms in (\ref{scheme:infinity}), we find balance properties equivalent to (\ref{eq:balance1})--(\ref{eq:balance2}), but with remainders $E_I$ and $F_I$ depending on $H_1$ and $H_2$ instead of $H_{i\geqslant \f{I}{2}}.$ One has to check \emph{a posteriori} that these last quantities are very small ; it is not the case in a standard calculation, but becomes true when the precision of the direct problem scheme increases. 

\subsection{Inverse problem: Quasi-Reversibility Approach}
\label{subsec:inverse}

In this section, we present a numerical scheme for the regularized inverse problem proposed in \cite{PZ}.
This problem leads to solving (\ref{eq:invfull}) taken at $y=2x$, that is
$$\left\{
\begin{array}{l}
\al \f{\p}{\p y} (B_{\e,\al} N_\e) + 4B_{\e,\al}(y) N_\e(y) = B_{\e,\al}\big(\f y 2\big) N_\e\big(\f y 2\big)+ \lb_{\e,\alpha} N_\ep\big(\f y 2\big)+2 \f{\p}{\p y}\biggl( N_\e\big(\f y 2\big)\biggr), \qquad y>0,
\\
\\
(B_{\e,\al} N_\e)(0)=0,
\end{array} \right.
$$
where $\alpha>0$ is the regularizing parameter and $\lb_{\ep,\alpha}$ is defined by (\ref{eq:conserv2:PZ}). This gives, in a discretized version, after dropping the index $\ep$,
\beq 
\lb_{\ep,\alpha}=\f{\sum N_i}{\sum x_i N_i + \f{\alpha}{4}\sum N_i}.
\label{def:lbep:Q}
\eeq

For the numerical discretization we set $H^Q_{-1}=0$ and also recall that $N_0=0$ and assume that the data satisfies $N_{I+1}=0.$ We use a standard upwind scheme for the differential term: 
\beq \label{eq:scheme:inverse}
\f{\alpha}{\Delta x}(H^Q_i-H^Q_{i-1})+4H_i^Q=H_{\f{i}{2}}^Q + L^Q_{\f{i}{2}}, 
\eeq
where we have defined the fractional indices as in the filtering approach  by (\ref{def:G1/2}), and here
$$L^Q_i=\lb_\ep N_i+\f{N_{i+1}-N_{i}}{\Delta x}.$$
If we neglect the terms $H^Q_{i\geqslant \f{I}{2}+1}$, we can easily verify 
a discrete version of the balance laws (\ref{eq:conserv1:PZ}) and (\ref{eq:conserv2:PZ}), equivalent to (\ref{eq:balance1})--(\ref{eq:balance2}).

\section{Numerical Tests}
\label{sec:result}

As input data, we take the values of the function $N$ obtained by the numerical solution of the direct problem in 
Section~\ref{subsec:direct}, we add a random noise uniformly distributed in $[-\f{\ep}{2},\f{\ep}{2}],$ and we enforce nonnegativity of the data
$$
N_\ep=\max(N+\ep \,r,0).
$$
We solve the direct problem on a regular grid of $I+1$ points, on an interval $[0,2L].$ We need $L$ large enough, such that it is possible to assume that $N(x\geqslant {L}) \approx 0$ and we have  checked it \emph{a posteriori.} Indeed, we have seen that this property is essential when we use  
the inverse schemes on a domain $[0,L]$ in order to verify the balance laws (\ref{eq:conserv1:PZ})--(\ref{eq:conserv2:PZ}). 
In other words, we solve the direct problem on a domain twice larger than for the inverse problem. In the numerical tests we take $L=4,$ and we show the numerical solution $N$ only on the interval $[0,L]$ since it is uniformly small on $[L,2L].$

We solve the inverse problem by the different methods on a regular grid of $I_1+1$ points on $[0,L],$  with $\Delta x_1={L}/{I_1}.$ This grid is taken ten times finer than the grid used for the direct problem, \emph{i.e.} we take $I_1=10 I.$ Since we have chosen $L$ large enough so that $N(x\geqslant {L})\approx 0,$ we have always obtained that indeed $H(x \geqslant {L})\approx 0.$

As before, we denote by $H^Q$ and $H^f$  the solution data $H$ obtained respectively by the quasi-reversibility method of Section~\ref{subsec:inverse} and by the first filtering approach (from zero) of Section~\ref{subsec:filter}. We also define a solution $H^{fQ}$ by mixing both methods, \emph{i.e.} by 
solving the following equation:
\beq
\label{eq:inverse:filter+inverse}
\left\{	\begin{array}{cl}
& \al \f{\p}{\p x} (B_{\e,\al} N_\e) (y) + 4B_{\ep,\alpha}(y)N_\ep(y) -B_{\ep,\alpha}(x)N_\ep(x)=\biggl(\f{\p}{\p x} N_\ep +\lb_{\ep,\alpha} N_\ep\biggr)* \rho_\alpha (x),\quad x\geqslant 0,
\\ \\
&B_{\ep,\alpha}(x=0)N_\ep(x=0)=0,
\end{array} \right.
\eeq
where $\lb_{\ep,\alpha}$ is defined by 
\beq\label{eq:conserv2:PZQ}
\lb_{\ep,\alpha}= \f{\int_0^{\infty} N_\ep * \rho_\alpha dx}{\int_0^{\infty} x N_\ep *\rho_\alpha dx + \f{\alpha}{4}\int_0^{\infty} N_\ep * \rho_\alpha dx}.
\eeq
The relative error is measured, as seen in Theorem~\ref{th:estim:filter} and in Theorem~5.1 of \cite{PZ}, by 
$$\delta^Q=\f{||B N_\ep - H^Q||^2_{l^2}}{||N_\ep||_{l^2}}, \qquad \delta^f=\f{||B N_\ep - H^f||^2_{l^2}}{||N_\ep||_{l^2}},\qquad \delta^{fQ}=\f{||B N_\ep - H^{fQ}||^2_{l^2}}{||N_\ep||_{l^2}} .
$$
We have divided by $||N_\ep||_{L^2}$ and not by $||N||_{H^2}$ because in practice we only know the entry data with noise.

In order to illustrate the accuracy of our method, we also compare it to a naive way (brute force) of considering the equation. Namely, we approximate $\f{\p} {\p x} N(x)$ by a second-order Euler scheme without regularization. It gives a solution $H^b$ by the same formula (\ref{eq:scheme:inverse}), where we simply take $\alpha=0.$

\

\paragraph{The Direct Problem.}

\begin{figure}[ht]
\begin{center}
\begin{minipage}{17cm}
\includegraphics[width=8cm, height=7cm]{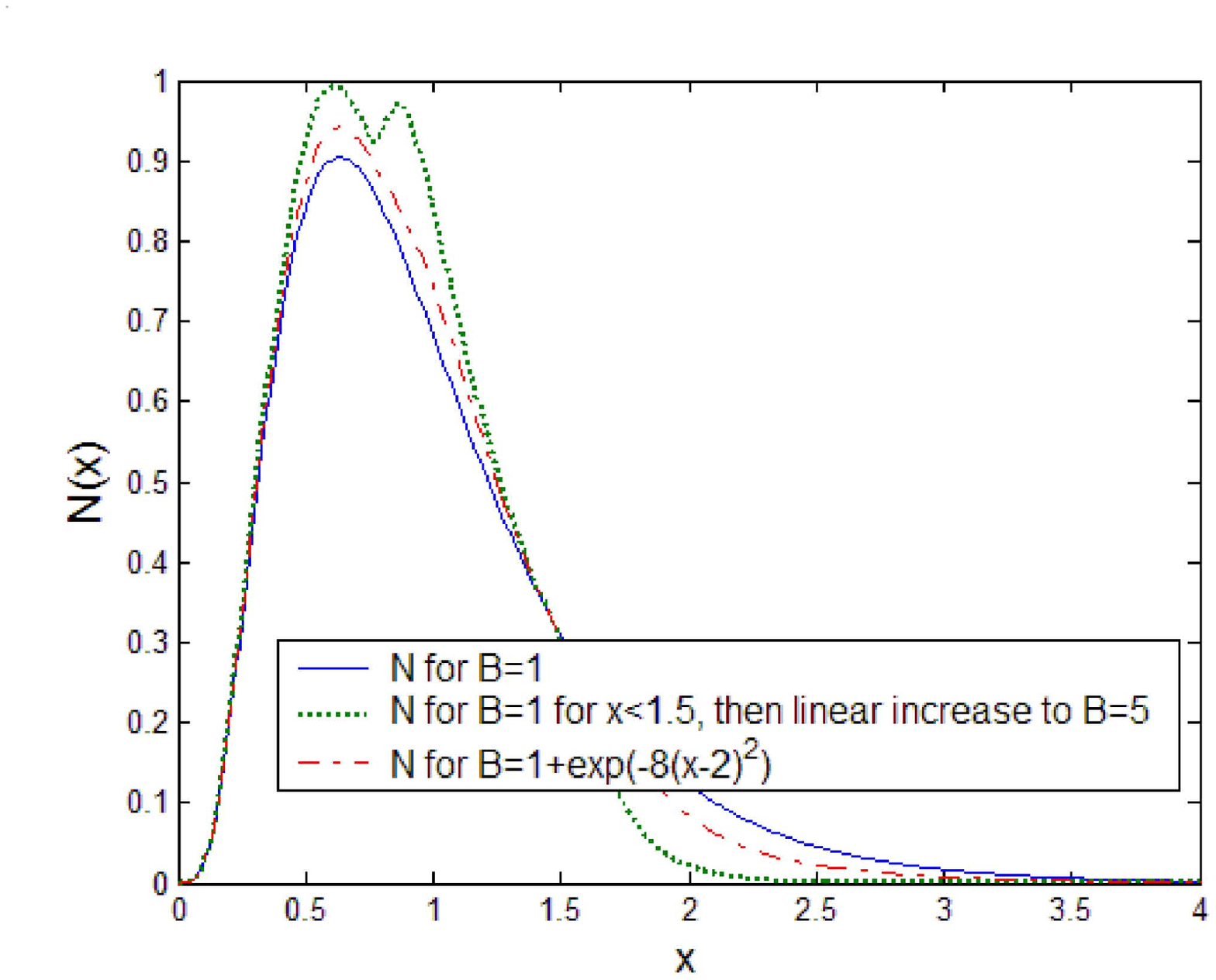} \quad \includegraphics[width=8cm,height=7cm]{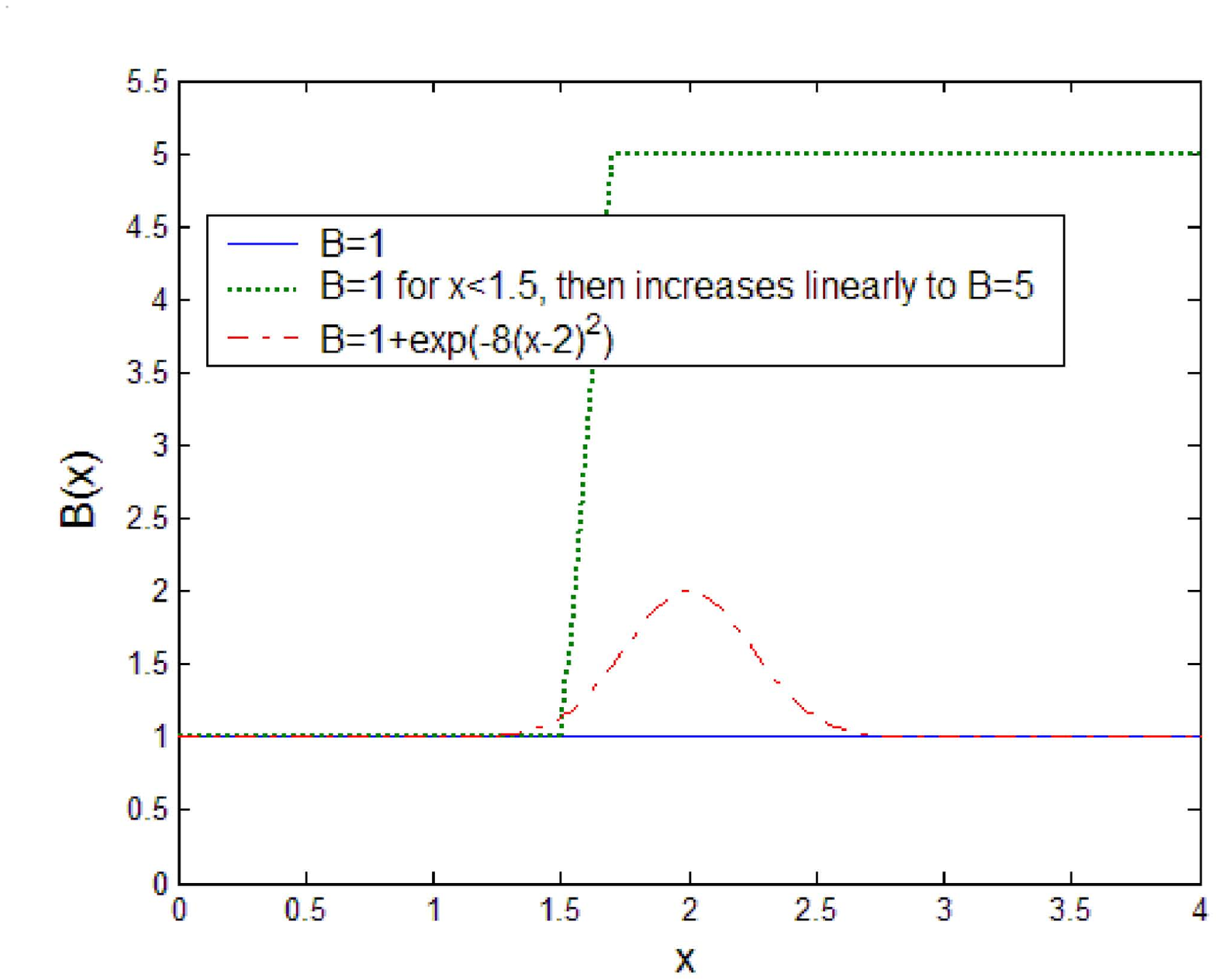}\end{minipage} \end{center}\vspace{-0.5cm}
\caption{\label{fig:ref:N} Solutions $N$ (left) obtained by the numerical resolution of Section \ref{subsec:direct} for the direct problem with three different division rates $B$ (right).}
\end{figure}

We have first tested the direct problem for various division rates $B$. Three different solutions $N$ for three given division rates $B$ are depicted in Figure \ref{fig:ref:N} with $800$ grid points. 

In the particular case when $B$ is constant, we can go further and evaluate the computational error. Then, we know that $\lb =B$ and the exact solution $N_{\rm exact}$  can be explicitly calculated, as shown in \cite{PR,BP}, by the formula:
\beq 
 N_{\rm exact}(x)=\bar{N}\sum\limits_{n=0}^\infty \alpha_n e^{-2^n B x}, 
\label{def:NM}
\eeq
where the coefficients are defined recursively by  $\alpha_0=1$ and $\alpha_n = (-1)^n \f{2\alpha_{n-1}}{2^n - 1}$, and $\bar N$ is chosen to ensure the mass one normalization. 
We take $B=1$ and obtain the continuous curve of Figure~\ref{fig:ref:N}.  We can measure here  the relative error by
$$\delta^D = \f{\|N-N_{\rm exact}\|_{l^1} }{\|N_{\rm exact} \|_{l^1}},
$$
where $N$ represents the numerical solution of Section~\ref{subsec:direct}. We choose this norm because for $B$ constant, the solution of the adjoint problem is $\phi=1$ and the General Relative Entropy Principle (\cite{MMP, BP}) gives us that this quantity decreases along the time iterations. Still for $800$ points, we obtain $\delta^{D}=7.7.10^{-3}.$  

\paragraph{The noiseless case ($\ep=0$).}

In the simplest case where the data is perfectly known, \emph{i.e.} for $\ep=0,$ we verify that the different schemes allow us to recover $B.$ Since the precision of the data is directly linked to the number of points used in the scheme, we run the codes with $1.000$ points for the direct problem  (below, we will take only $100$ points).

\begin{figure}[ht]
\begin{minipage}[b]{0.46\linewidth}
\includegraphics [width=\linewidth, height=7cm]{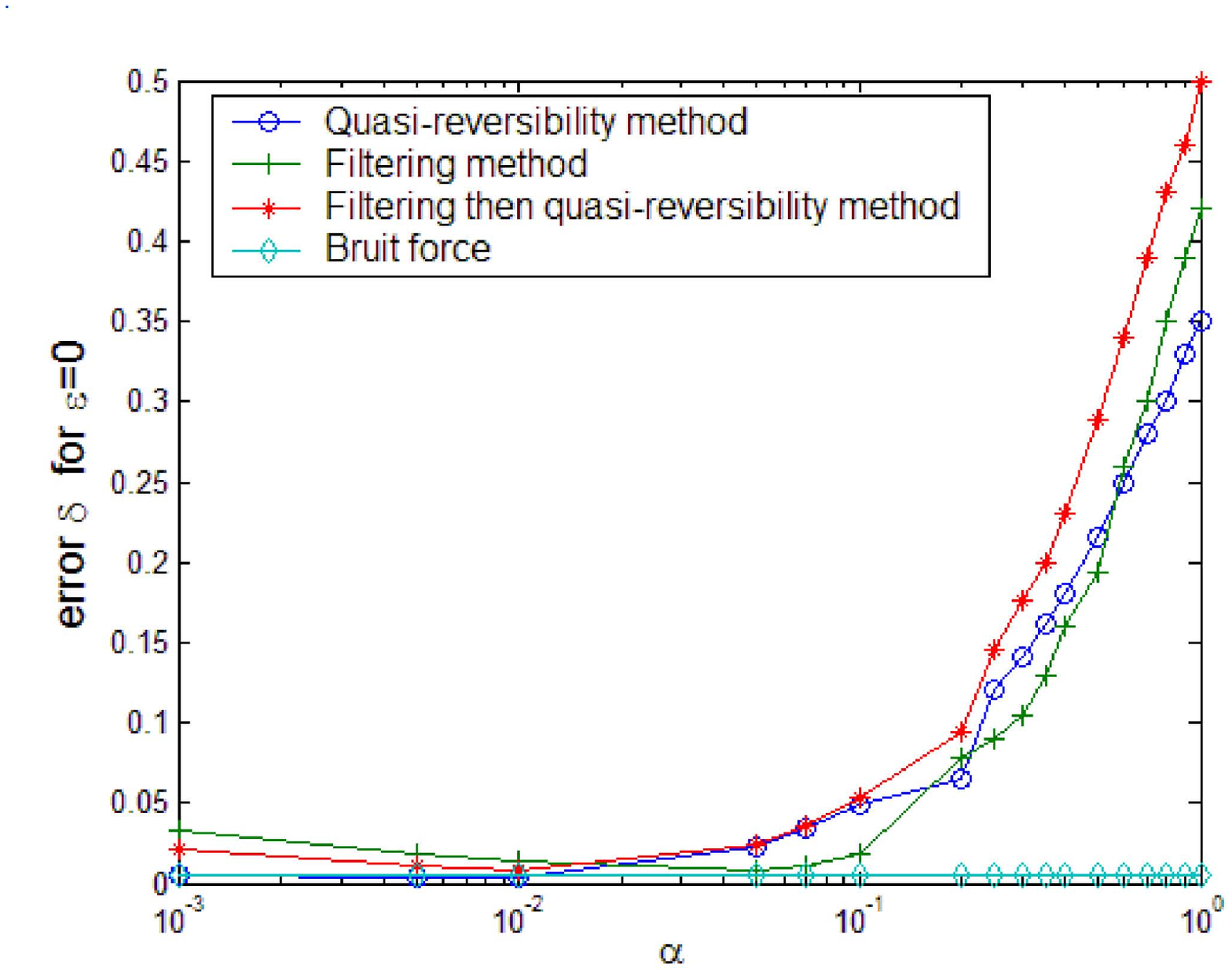}
\caption{\label{fig:table:ep0} For $\ep=0$, numerical errors obtained with the different methods for the inverse problem.}
\end{minipage} \hfill
\begin{minipage}[b]{0.46\linewidth}
\includegraphics [width=\linewidth, height=7cm]{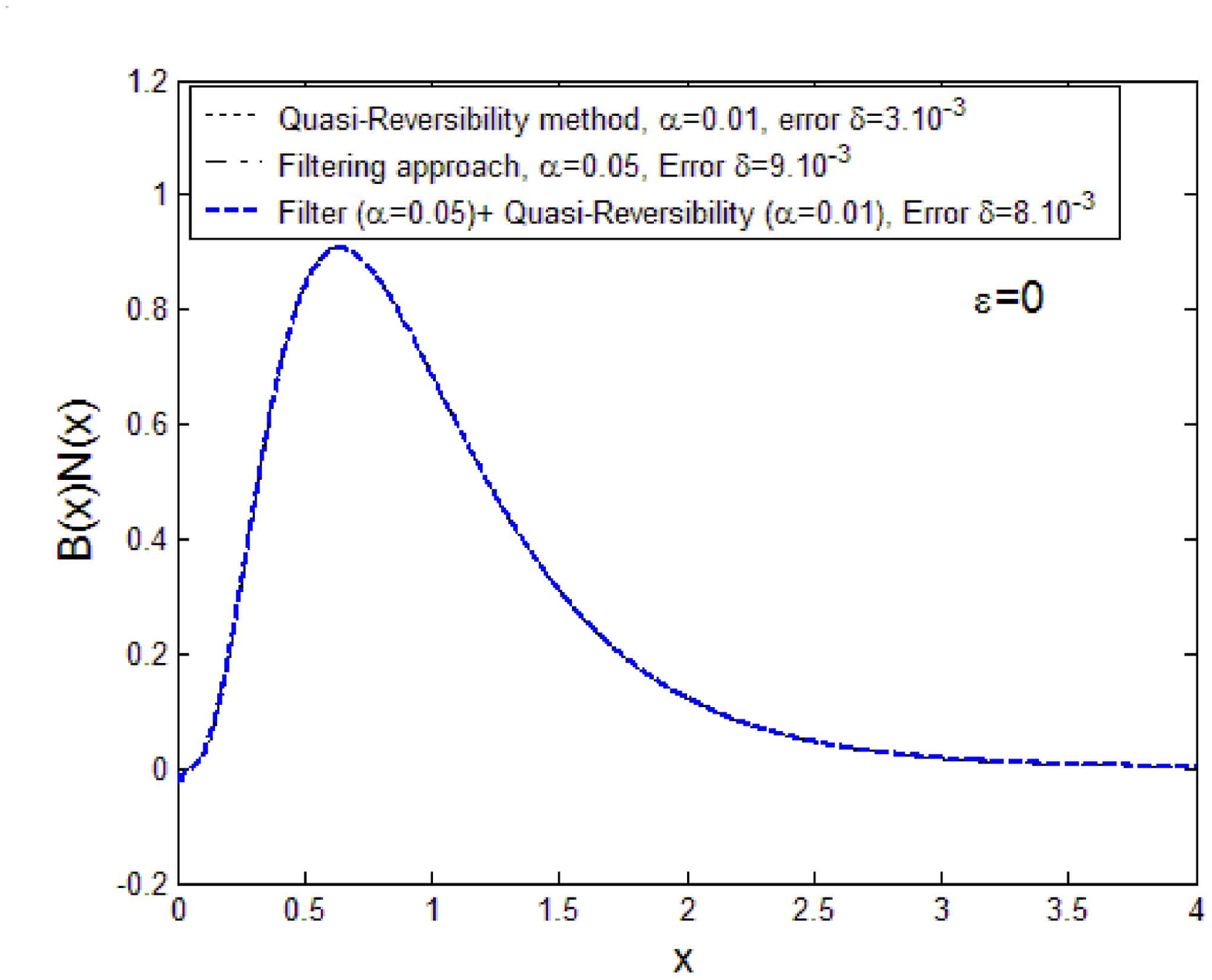}
\caption{\label{fig:ep0:NB1} Numerical reconstruction of $B.N$ obtained by each method for the inverse problem when $B=1,$ $\ep=0$ and $\alpha=0.01.$}
\end{minipage} \hfill
\end{figure}

\begin{figure}[ht]
\begin{minipage}[b]{0.46\linewidth}
\includegraphics [width=\linewidth, height=7cm]{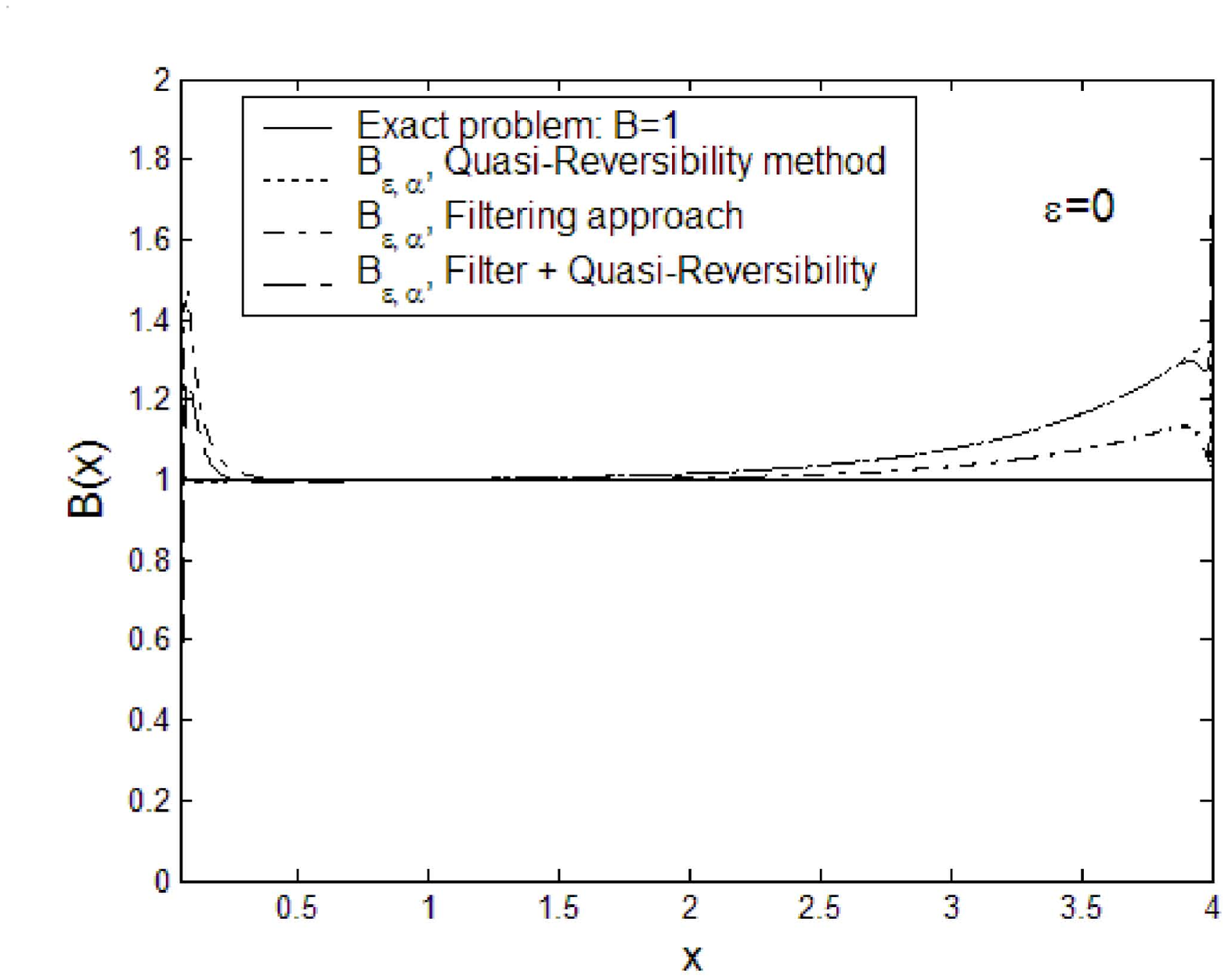}
\caption{\label{fig:ep0:B1} Reconstructed division rate $B$ using the three inverse methods, for $\ep=0$, $\al=0.01$ with $N$ computed from $B=1.$}
\end{minipage} \hfill
\begin{minipage}[b]{0.46\linewidth}
\includegraphics [width=\linewidth, height=7cm]{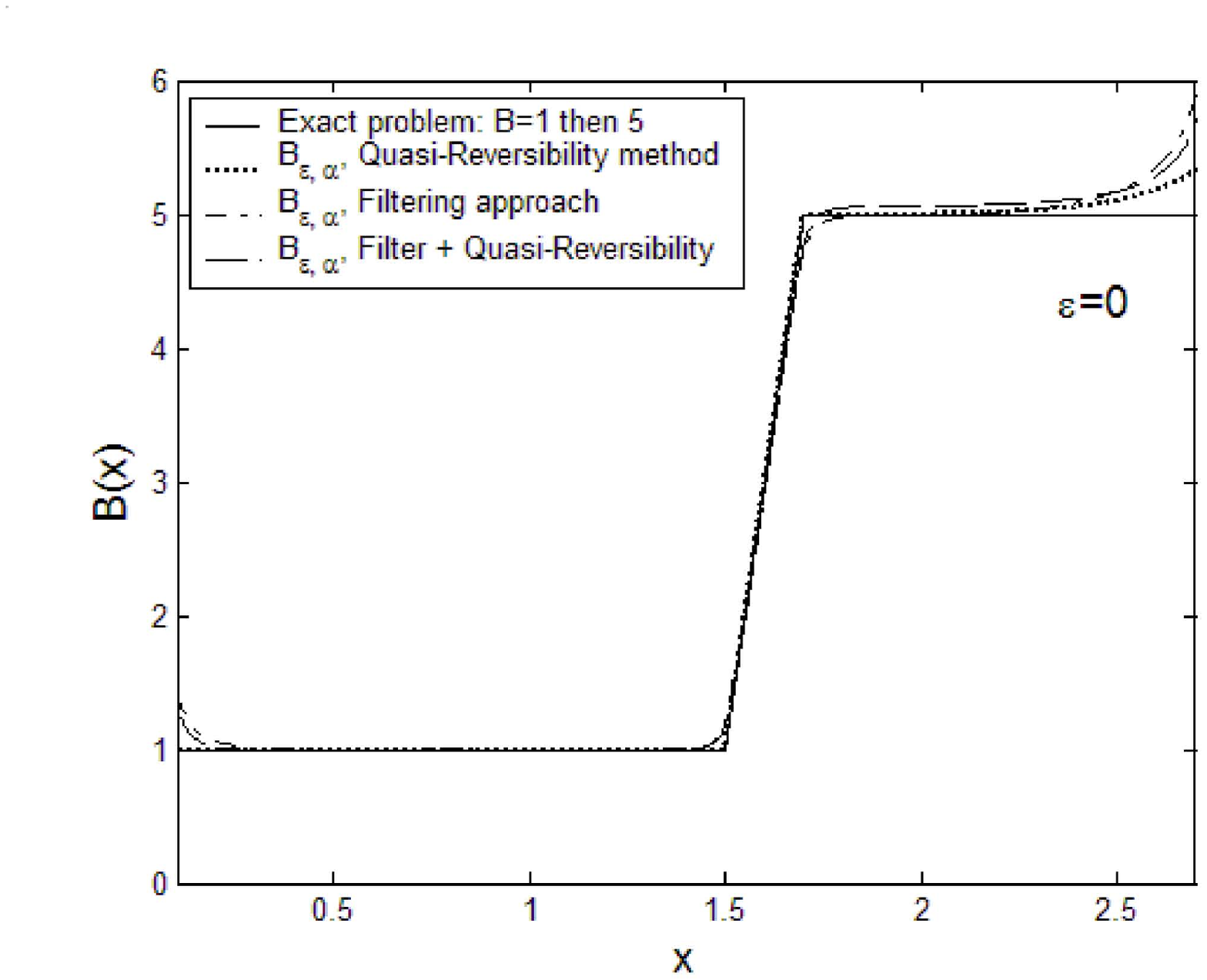}
\caption{\label{fig:ep0:Bexp} Reconstructed division rate $B$, for $\ep=0,$ $\al=0.01$ and a jump $B=1$ to $5$ as in  Figure \ref{fig:ref:N}}
\end{minipage} \hfill
\end{figure}

\begin{figure}[ht]
\begin{minipage}[b]{0.46\linewidth}
\includegraphics [width=\linewidth, height=7cm]{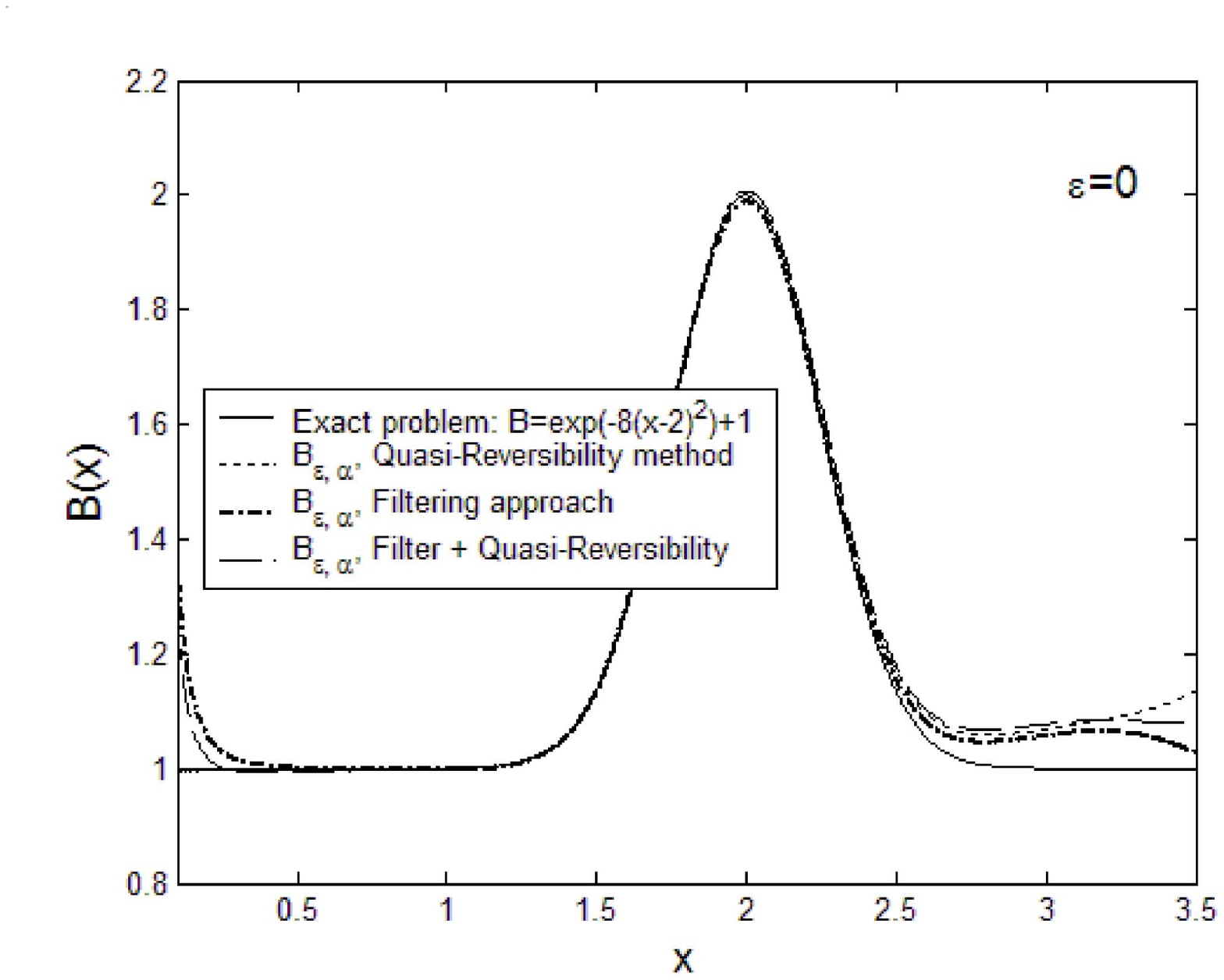}
\caption{\label{fig:ep0:Bcos} Reconstructed division rate $B$, for $\ep=0,$ $\al=0.01$ and $B=1+exp(-8(x-2)^2)$.}
\end{minipage} \hfill
\begin{minipage}[b]{0.46\linewidth}
\includegraphics [width=\linewidth, height=7cm]{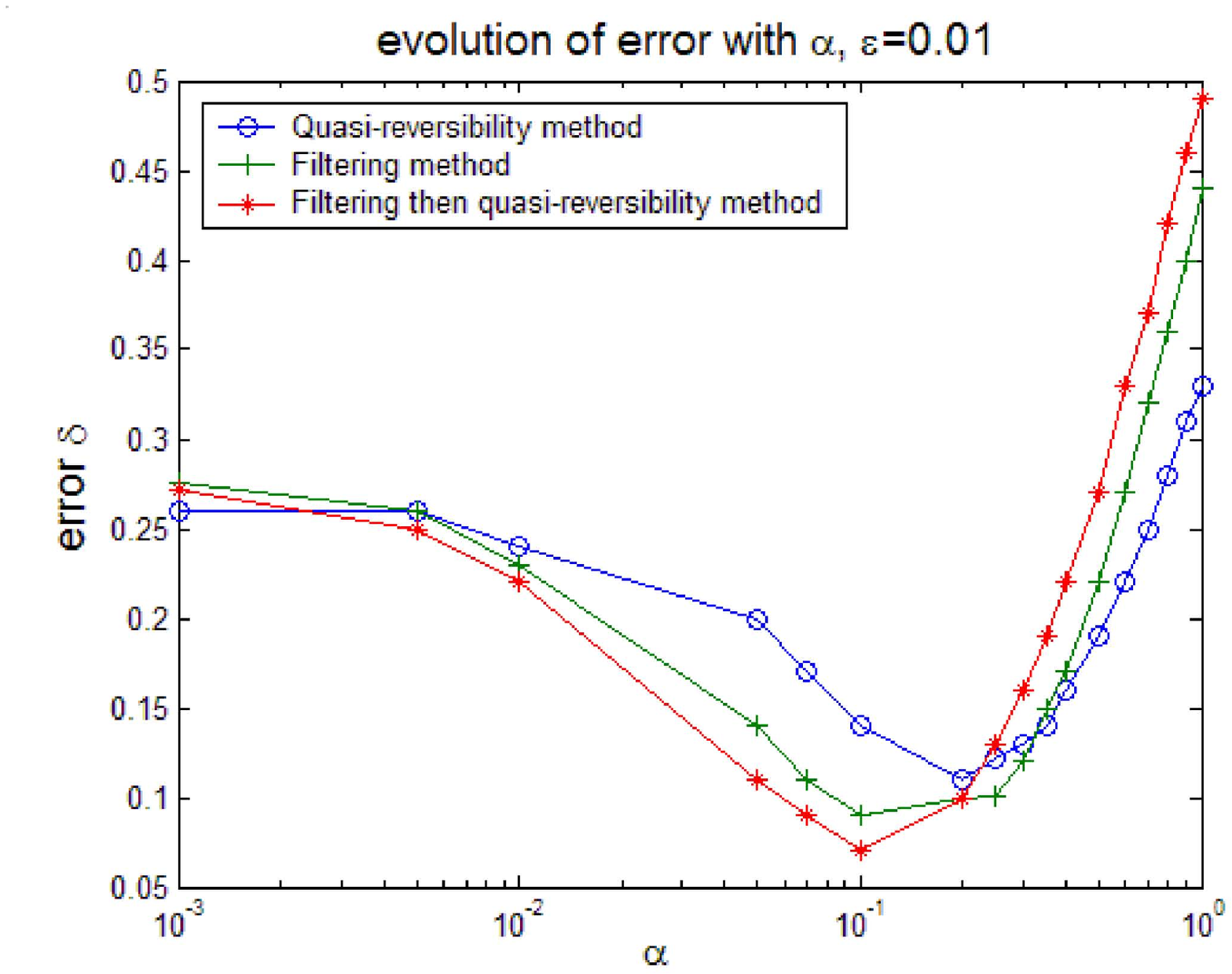}
\caption{\label{fig:table:ep0_01} Numerical error when $\ep=0.01$ for the different methods. }
\end{minipage} \hfill
\end{figure}

We test several values of $\alpha$ and we use the three functions $B$ of Figure \ref{fig:ref:N} for each method for the inverse problem. The error estimate is found to depend on the method used but not significantly on the division rate $B$. Therefore  
we have drawn in Figure~\ref{fig:table:ep0} the  average error estimates for the three division rates $B$. In Figure~\ref{fig:ep0:NB1} we have depicted the products $B.N$ in the case $B=1$ and $\alpha=0.01$ (other cases are similar): it shows that the precision obtained is satisfactory. In Figures~\ref{fig:ep0:B1}, \ref{fig:ep0:Bexp} and \ref{fig:ep0:Bcos} we have drawn the approximations of B in each of the three cases, calculated only for $N>0.01$ (indeed, for $N$ too small the division leads to insignificant results on $B$). 

Not surprisingly, the brute force method reveals to be satisfactory, with an error estimate of $\delta^b=1.3.10^{-2},$ since we are in the case where $N$ is very regular. The filtering method can reach this level of error for $\alpha=10^{-2}$ but cannot go further. 
However, both the quasi reversibility method and the mixed method given by Equation (\ref{eq:inverse:filter+inverse}) improve it with minimum values $\delta^Q = 6.9.10^{-3}$ and $\delta^{fQ}=6.5.10^{-3}$ reached for $\alpha=10^{-2}.$

\paragraph{Link between the noise level $\ep$ and the regularization parameter $\alpha$.}

\begin{figure}[ht]
\begin{minipage}[b]{0.46\linewidth}
\includegraphics [width=\linewidth, height=7cm]{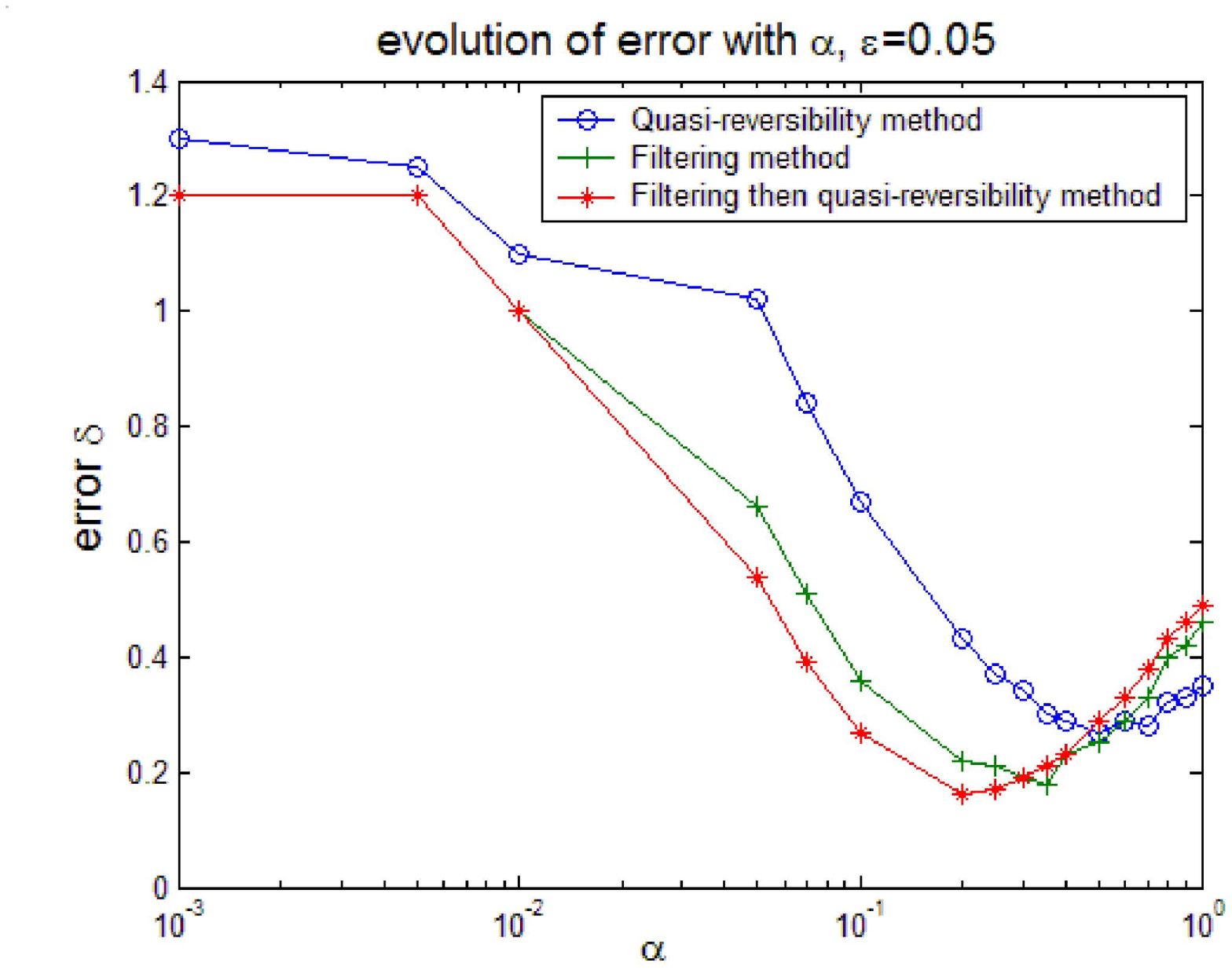}
\caption{\label{fig:table:ep0_05} Numerical error when $\ep=0.05$.}
\end{minipage} \hfill
\begin{minipage}[b]{0.46\linewidth}
\includegraphics [width=\linewidth, height=7cm]{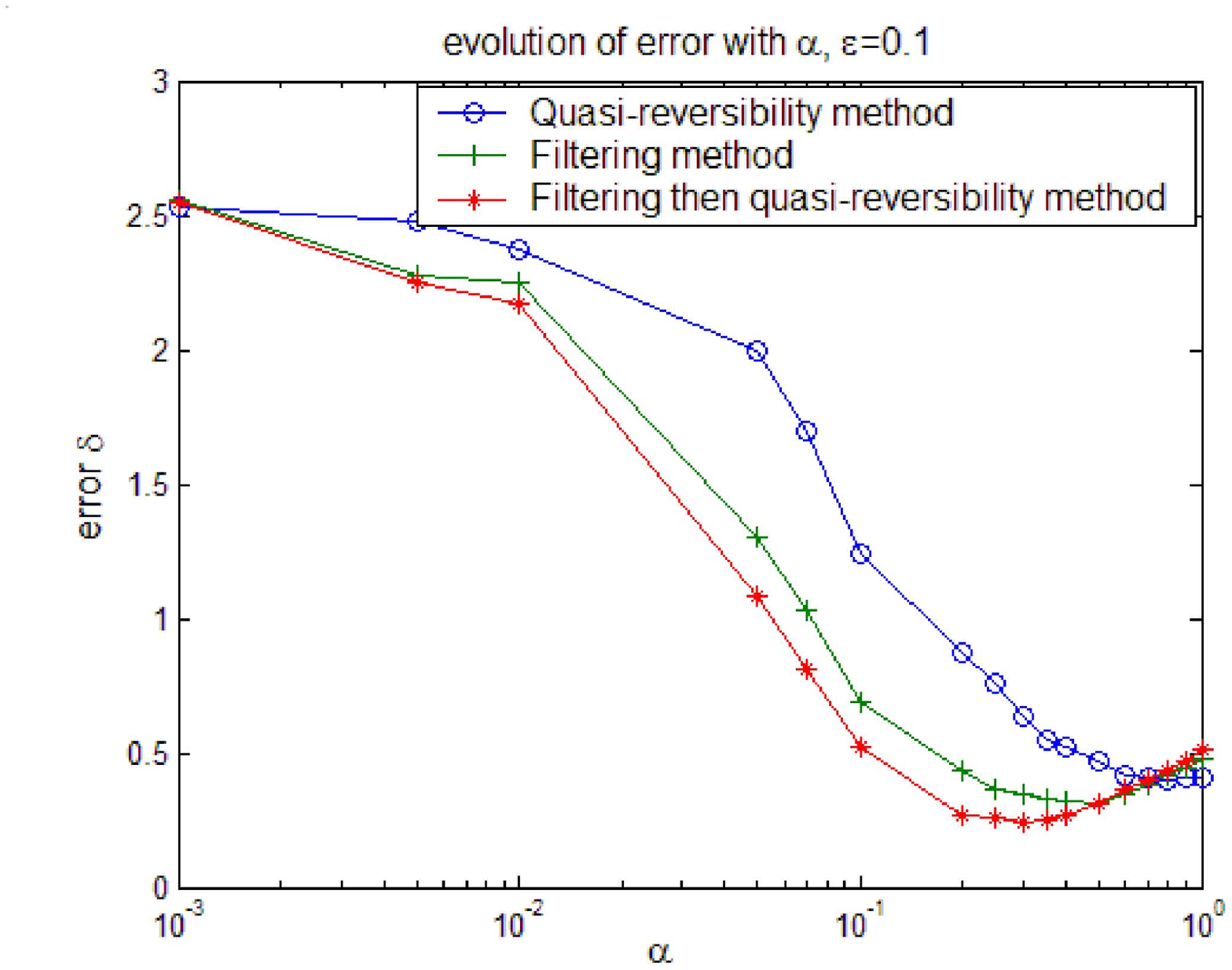}
\caption{\label{fig:table:ep0_1} Numerical error when  $\ep=0.1$.}
\end{minipage} \hfill
\end{figure}

\begin{figure}[ht]
\begin{center}
\begin{minipage}{17cm}
\includegraphics [width=8cm, height=7cm]{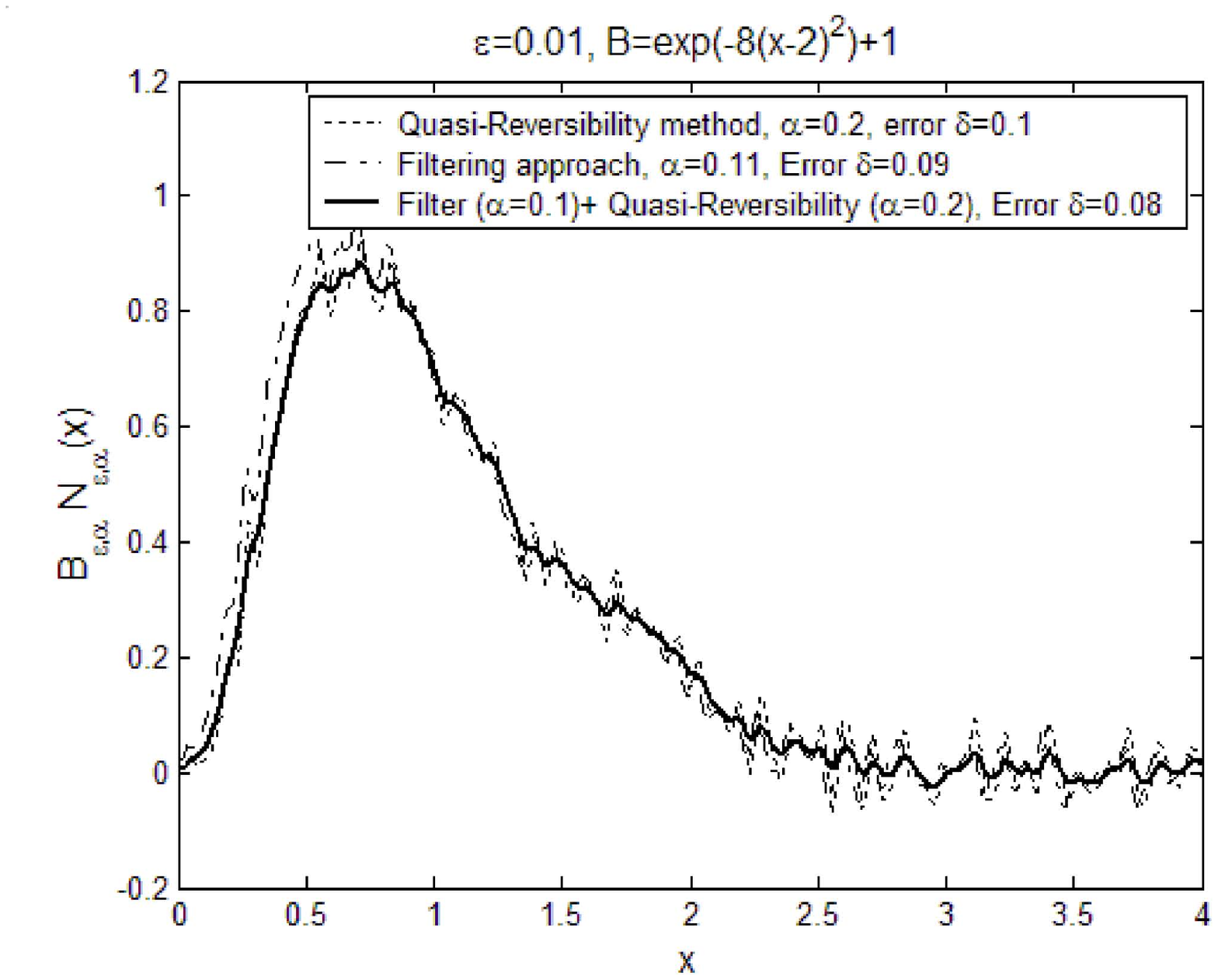} \quad 
\includegraphics [width=8cm, height=7cm]{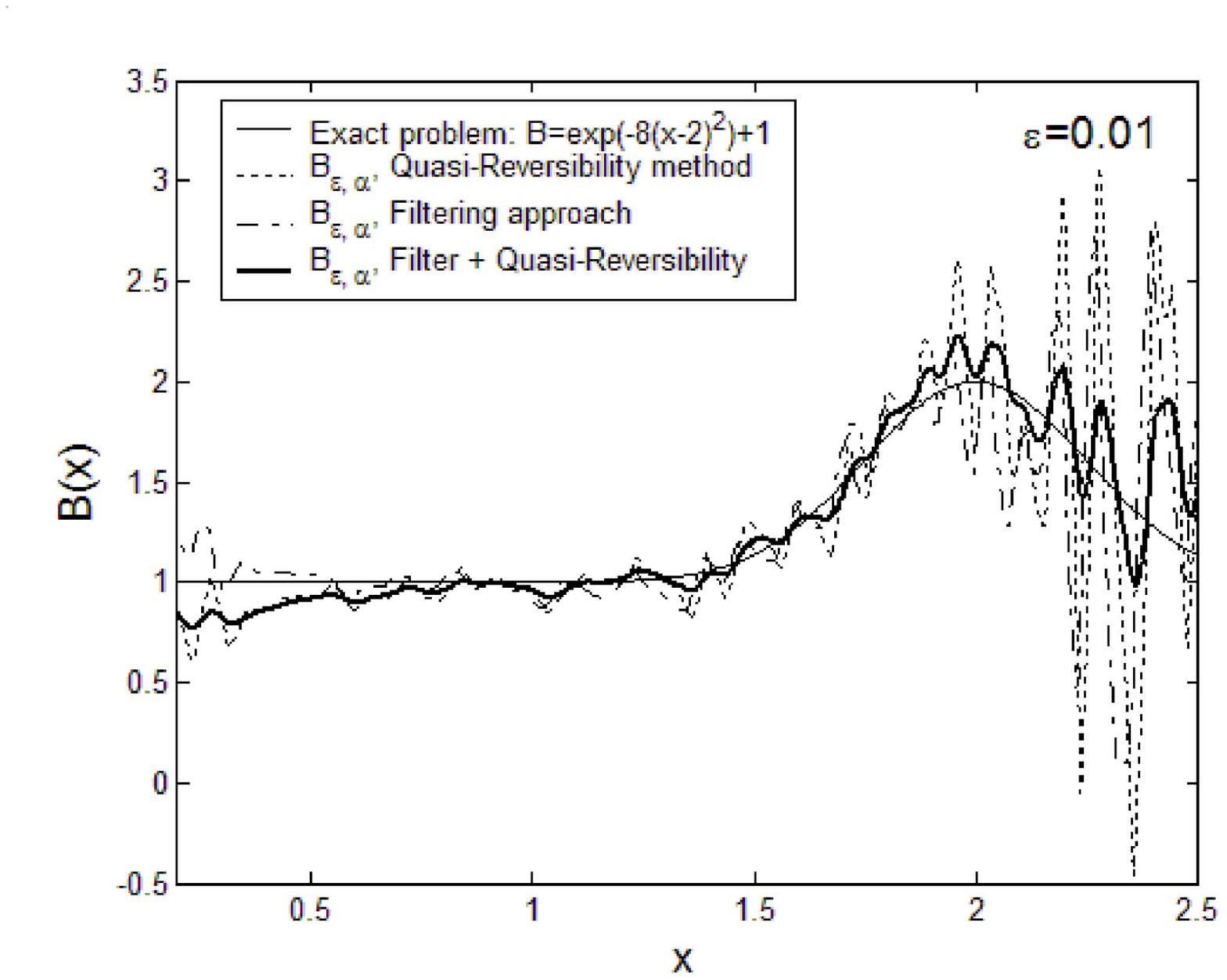}\end{minipage} \end{center}\vspace{-0.5cm}
\caption{\label{fig:ep0_01:B} In the case  $\ep=0.01,$ $\alpha=0.05$, $B=1+e^{-8(x-2)^2}$, numerical solution $B.N$ (left) and $B$ (right) by the different methods.}
\end{figure}

\begin{figure}[ht]
\begin{minipage}[b]{0.46\linewidth}
\includegraphics [width=\linewidth, height=7cm]{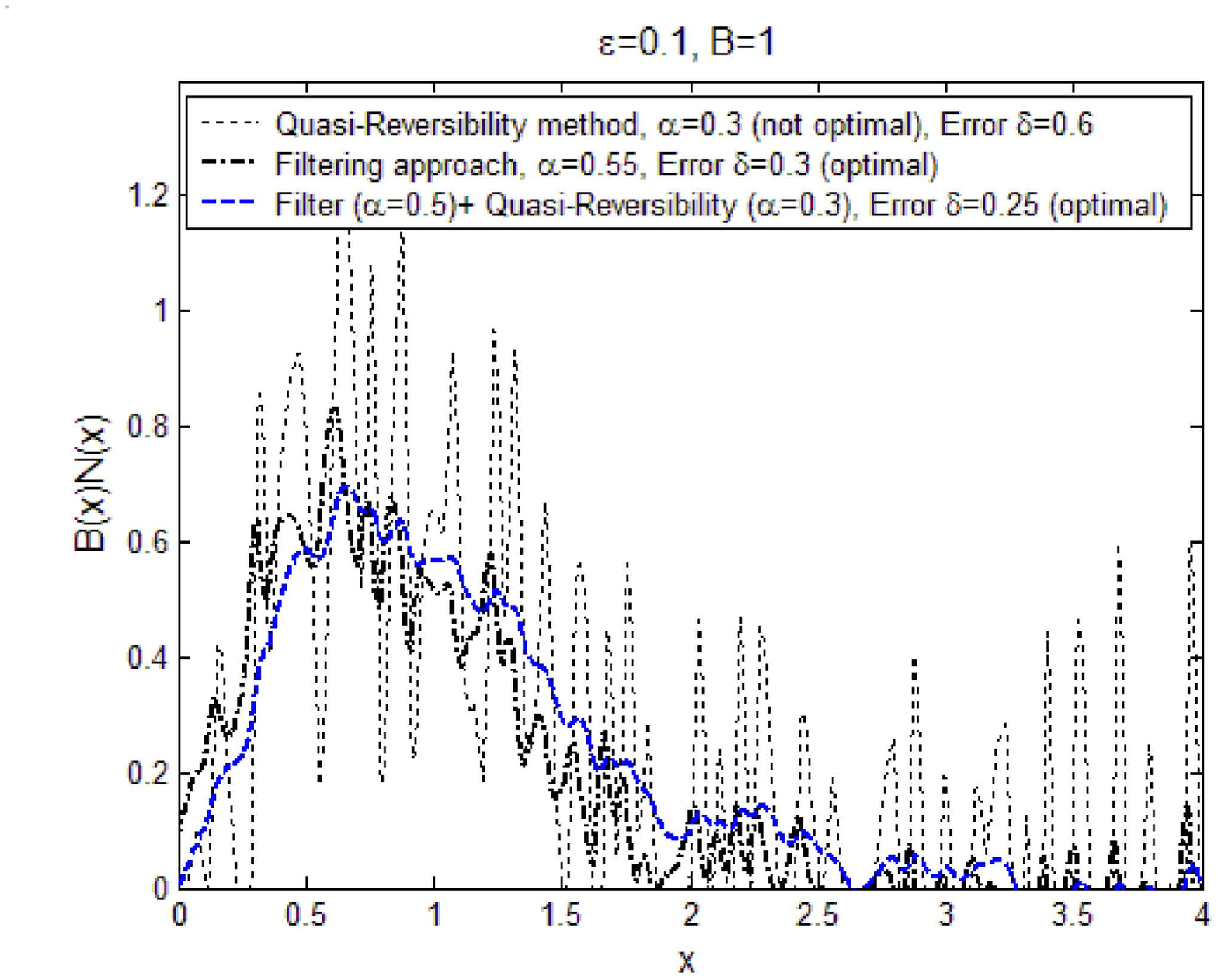}
\caption{\label{fig:ep0_1:alp0_3BN} In the case $\ep=0.1$ and $B=1$, the numerical solution $B.N$ by the different methods. \qquad  \qquad \qquad \qquad \qquad \qquad \qquad \qquad \qquad \qquad \qquad }
\end{minipage} \hfill
\begin{minipage}[b]{0.46\linewidth}
\includegraphics [width=\linewidth, height=7cm]{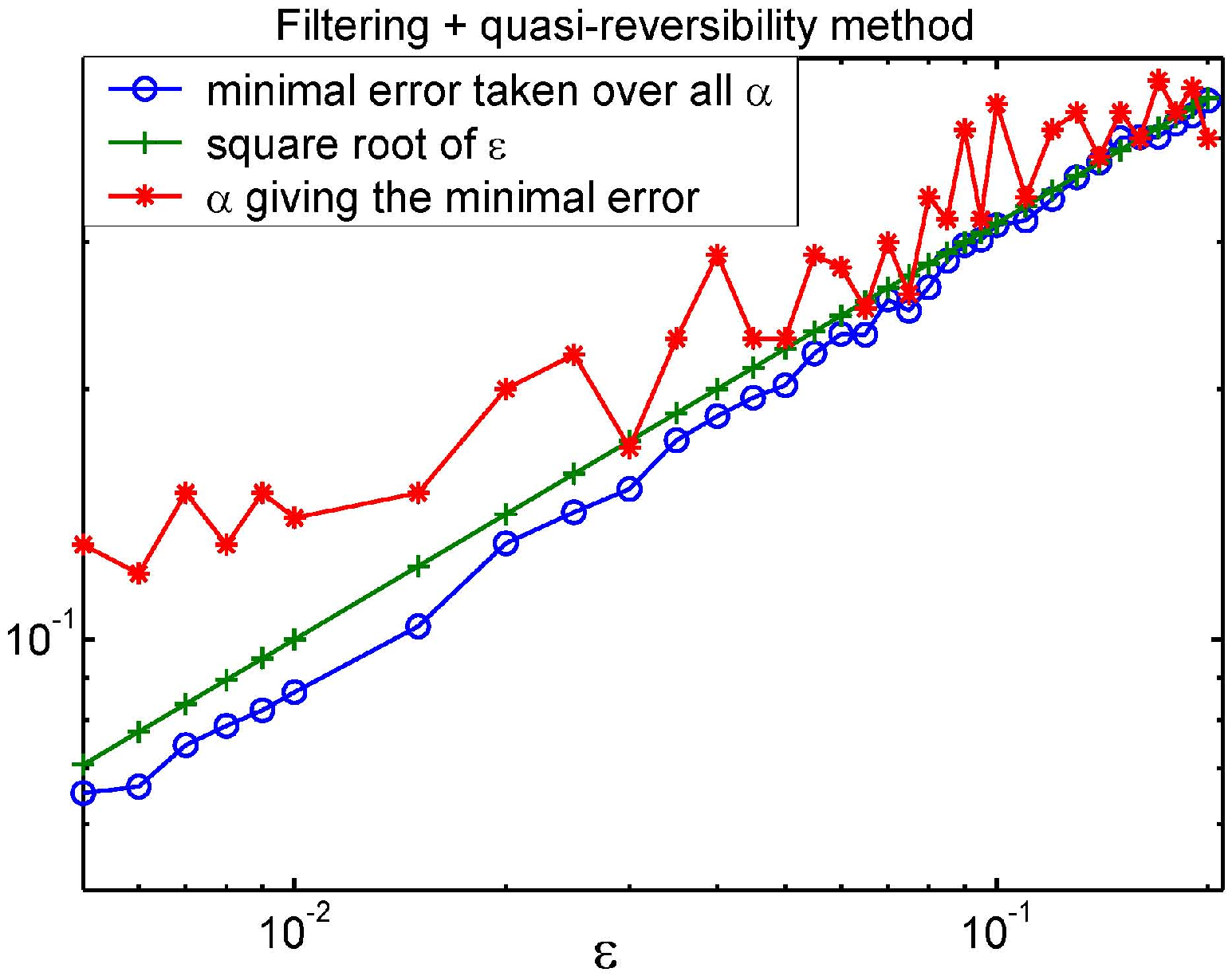}
\caption{\label{fig:QRF} Filtering and quasi-reversibility method: minimal error as a function of the noise level $\ep$ for the optimal value  $\alpha$, with a  comparison to the theoretical curve $\sqrt{\e}$.}
\end{minipage} \hfill
\end{figure}

For noise levels $\ep=0.01,$ $\ep=0.05$ and $\ep=0.1$ respectively, the Figures~\ref{fig:table:ep0_01}, ~\ref{fig:table:ep0_05} and~\ref{fig:table:ep0_1} give the curves  $\e$ as a function of $\alpha$ for the three inverse methods.  We compare the reconstructed division rates $B$ in  Figures~\ref{fig:ep0_01:B} and \ref{fig:ep0_1:alp0_3BN}.

\begin{figure}[ht]
\begin{center}
\begin{minipage}{17cm}
\includegraphics [width=8cm, height=7cm]{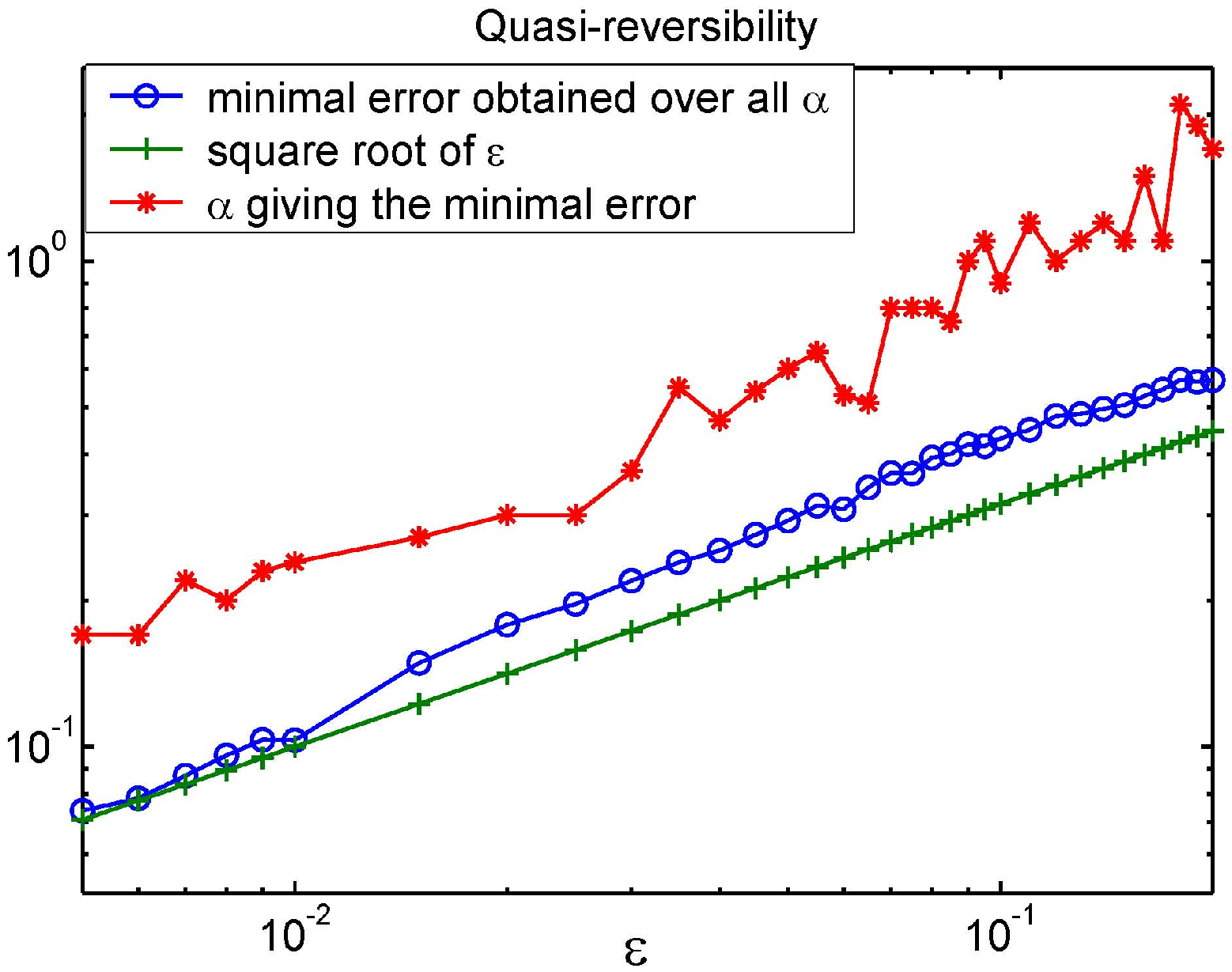}
\includegraphics [width=8cm, height=7cm]{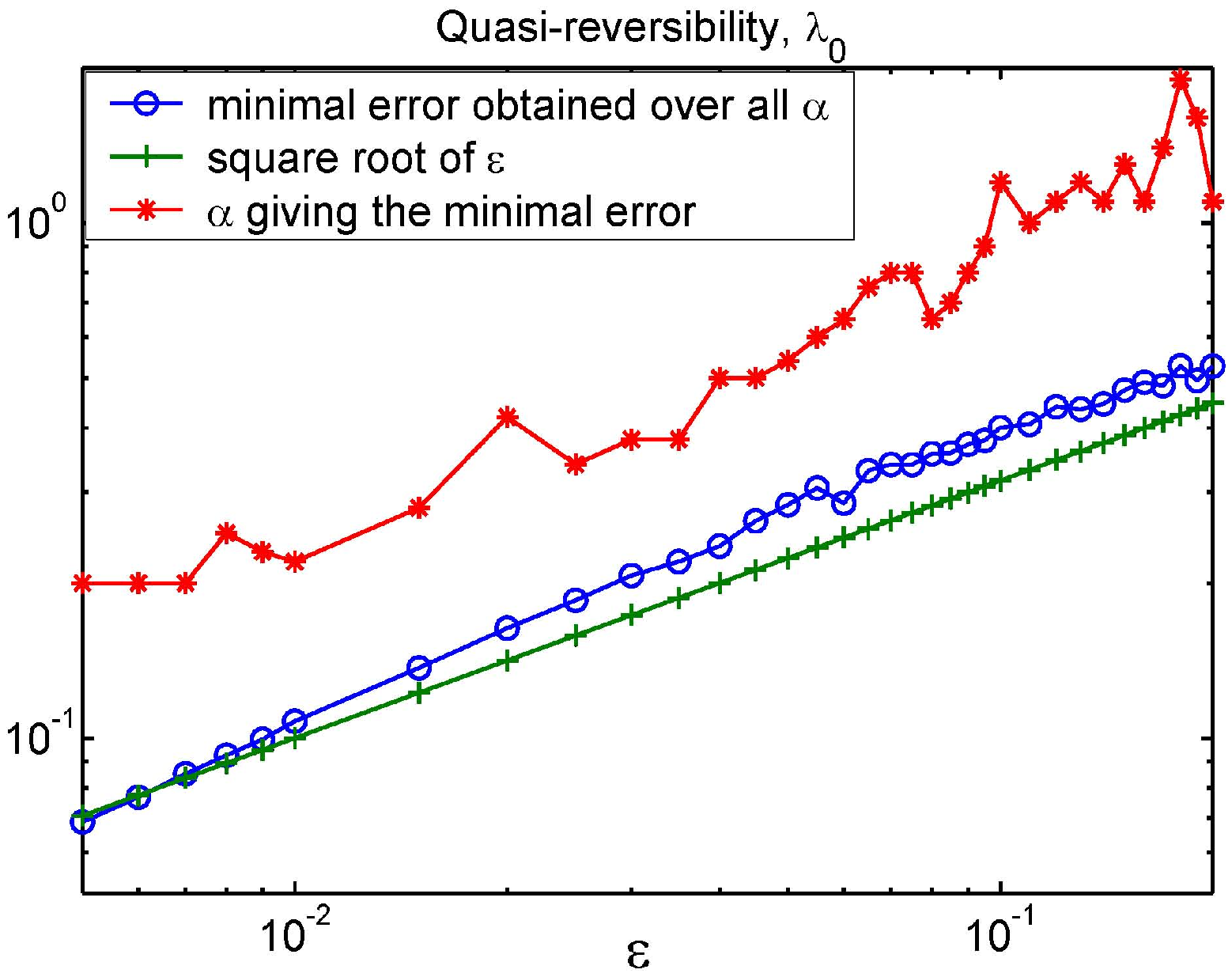}
\end{minipage} 
\end{center} \vspace{-0.5cm}
\caption{ \label{fig:QR} Quasi-reversibility method, with (left) $\lb_{\ep,\alpha}$ given by relation (\ref{eq:conserv2:PZ}) or (right) \emph{a priori} knowledge of $\lb_0:$  minimal error and optimal regularization parameter $\alpha$ as functions of the noise level $\ep$, with a  comparison to the theoretical curve $\sqrt{\e}$. We see that the \emph{a priori} knowledge of $\lb_0$ does not improve the speed of convergence of the scheme.}
\end{figure}

\begin{figure}[ht]
\begin{minipage}[b]{0.46\linewidth}
\includegraphics [width=\linewidth, height=7cm]{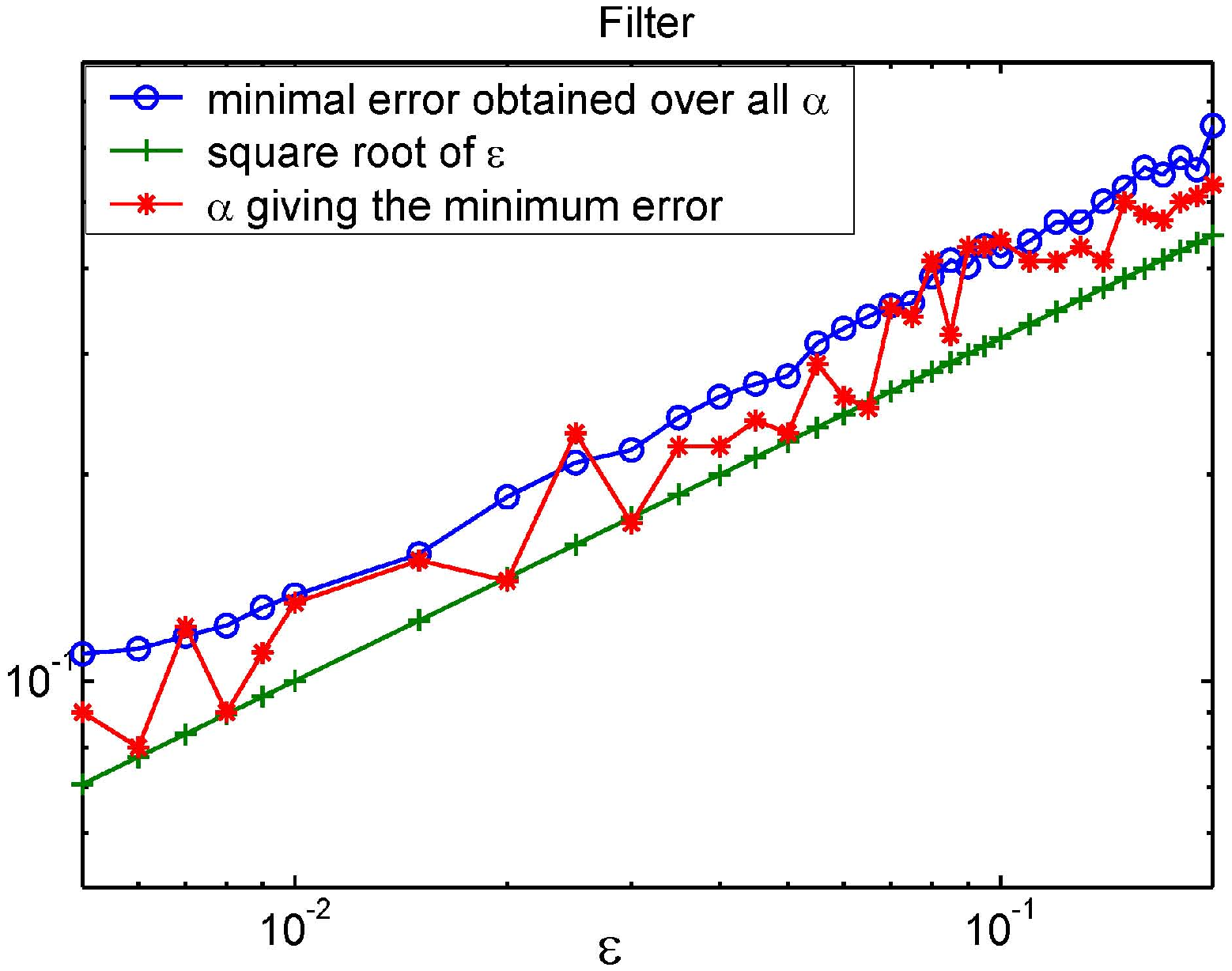}
\caption{\label{fig:filter} Filtering method for standard levels of noise.}
\end{minipage} \hfill
\begin{minipage}[b]{0.46\linewidth}
\includegraphics [width=\linewidth, height=7cm]{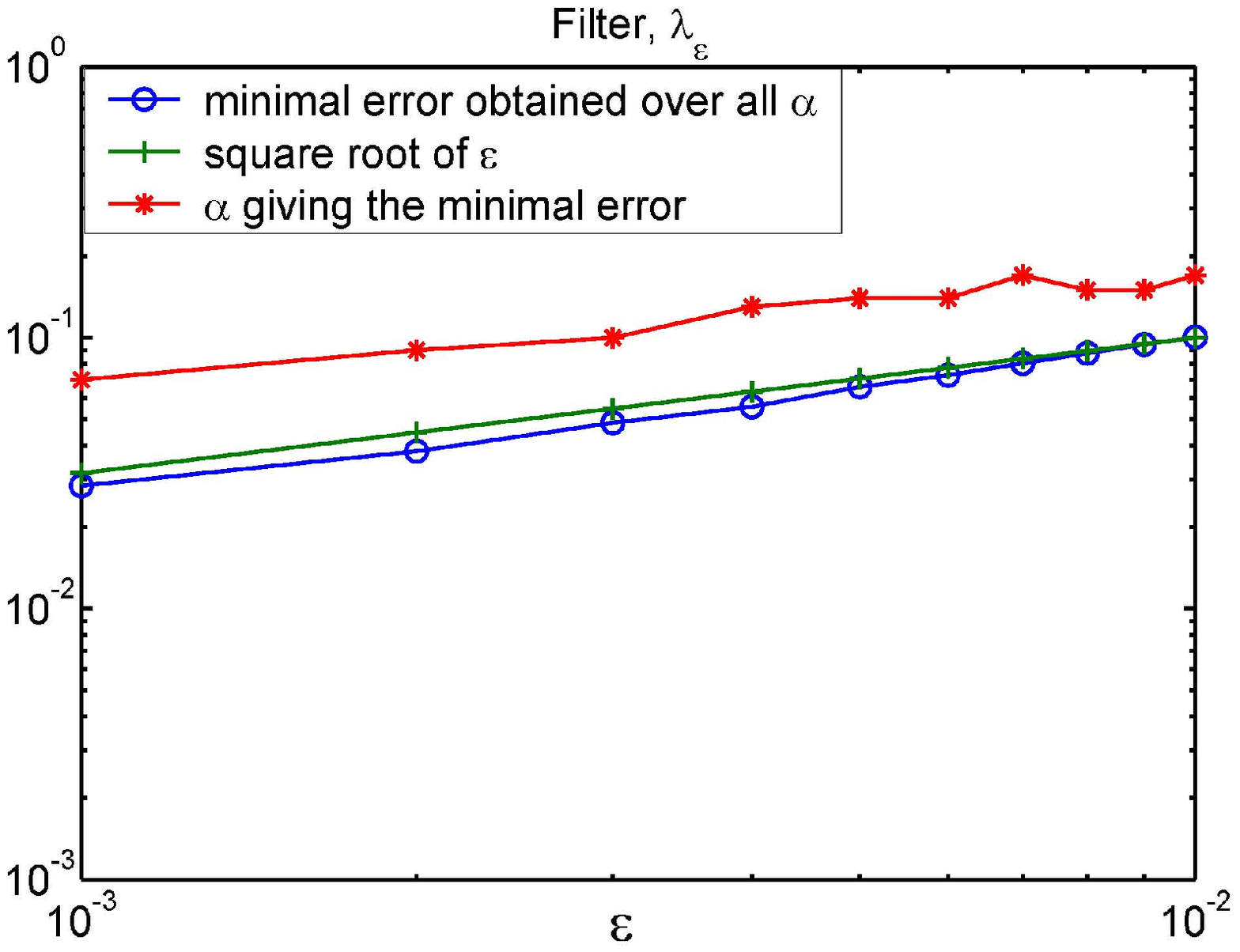}
\caption{\label{fig:filter:small} Filtering method for smaller values of $\ep$ and increased number of points.}
\end{minipage} \hfill
\end{figure}

Each of the  error curves  presents a minimum for  an optimal value of $\alpha$, as expressed by estimate (\ref{ineq:estim:filter}) for instance. 
In Figures \ref{fig:QRF}, \ref{fig:QR}, \ref{fig:filter} and \ref{fig:filter:small}, we have compared three curves, drawn in a log-log scale: $\sqrt{\ep}$ to serve as a reference curve, $f(\ep)=\min\limits_{\alpha} \delta (\alpha,\ep),$ and $g(\ep)=\text{arg} \min\limits_\alpha \delta(\alpha,\ep).$ One can see that for each method, these three curves have comparable slopes ($\f{1}{2}$ on a log-log scale): they show that even though the combination of filtering and quasi-reversibility method improves the optimal errors in absolute value, it does not change the order of convergence of the approximation, which remains of order $O(\sqrt{\ep}).$ Figure \ref{fig:filter:small} gives also the convergence of the filtering method for much smaller values of $\ep$ (for which an increased number of $500$ points has been taken, in order to avoid numerical bias): the comparison with $\sqrt{\epsilon}$ is there particularly evident, and we have obtained similar curves for the two other methods. 
The speed of convergence is though in complete accordance with the theoretical estimate (\ref{ineq:estim:filter}).

\paragraph{Influence of the choice of $\lb_0$ instead of $\lb_{\ep,\alpha}.$}

To evaluate the influence of the error term due to the distance $|\lb_{\ep,\alpha} -\lb_0|,$ we compare the curves obtained respectively by taking on the inverse code the exact $\lb_0$ or the value $\lb_{\ep,\alpha}$ expressed by the balance laws. They are drawn in Figure \ref{fig:QR} for the quasi-reversibility method. They show that even though the \emph{a priori} knowledge of $\lb_0$ improves the error in absolute value, it does not change the order of convergence of the scheme. Thus it is in complete accordance  with Estimate (\ref{ineq:estim:filter}).

\section{Conclusions}
We have considered size-structured equations connected to several areas of biology from cell division to prion proliferation by aggregation and fragmentation. We have addressed the numerical efficiency of some inverse problem solution methods to tackle the problem of  recovering the division rate from the size distribution of cells. The latter involves
a dilation equation with a singular right-hand side that needs regularization for actual implementation. For that purpose, we have introduced a filtering method and proved its convergence for noisy data. This method brings in an operator that has a non-trivial kernel and we have selected a numerical approximation that is able to recover the natural solution we want to reach. 

The implementation of the inverse algorithm, based on the filtering method, confirms the convergence analysis. In particular, there is an optimal regularization parameter as can be seen in the 
graphs of Figures~\ref{fig:table:ep0_01}, \ref{fig:table:ep0_05} or~\ref{fig:table:ep0_1} for instance. 
Comparison with a quasi-reversibility method introduced earlier leads to the conclusion that a combination of filtering  and quasi-reversibility methods seems to be more efficient because the oscillations are reduced, 
but without improving the rate of convergence. 

We also analyzed the impact of using the exact value of $\lb_0$ or the $\lb_\ep$ on the different solutions of the inverse problem.  In our simulations, the difference between using $\lb_0$ or $\lb_\ep$ seemed to be immaterial as far as the accuracy of method is concerned.
This is in perfect accordance with the theoretical estimate (\ref{ineq:estim:filter}). 

%
%
%
%
%
%

The above remarks open several directions for continuation and extension of the present work.  
On the practical side, the present work sets the stage for the use of experimental data
either from the existing literature or from more recent biological experiments.
On the theoretical side, the possibility of improving the convergence by combining the 
filtering and quasi-reversibility methods should be investigated further. 

Finally, we point out that although the 
Tikhonov method is  more standard,  we did not study it so far because it seems more time consuming. Indeed, iterations are needed to solve both the direct problem and the inverse one. To overcome such difficulty a completely new theory has to be developed so as to suit the particular structure of our model.  This provides yet another direction for future work. 

\section*{Acknowledgments}
The authors were supported by the CNPq-INRIA agreement INVEBIO. 
JPZ was supported by CNPq under grants 302161/2003-1 and
474085/2003-1. JPZ and BP are thankful to the RICAM special semester
and to the International Cooperation Agreement Brazil-France.  Part of this work was conducted 
during the Special Semester on Quantitative Biology Analyzed by Mathematical Methods, October 1st, 2007
-January 27th, 2008, organized by RICAM, Austrian Academy of Sciences.

\appendix

\section{Well-Posedness of Functional Equation Associated to the Inverse Problem.}

We have seen that the regularization method for the inverse problem relies mostly on solving the equation
\beq 
4 H(2x) - H(x) = L(x), \qquad x \geq 0.
\label{eq:frag}
\eeq 
Eventhough this equation is formally very simple, its analysis reveals some complexity. It may admit several solutions in general. Among them, we can mention two with simple representation formulas (we leave to the reader to check they are indeed formally solutions)
\beq\label{form:H1}
H^{(1)}(x)=  \sum\limits_{n=1}^{+\infty} 2^{-2n} L(2^{-n} x) .
\eeq
\beq \label{form:H0}
H^{(2)} (x) = -\sum\limits_{n=0}^{+\infty} 2^{2n} L(2^n x),
\eeq
To clarify this issue and motivate our choice of a solution, we first state general results concerning solutions to \fer{eq:frag} and then come back to our original problem (\ref{eq:exact:inverse}).

We first mention the following
\begin{proposition}\label{prop:pbmH}
Let $L\in L^2(\R_+, x^p dx),$ with $p\neq 3$, then there exists a unique solution $H\in L^2(\R_+, x^p dx)$ to \fer{eq:frag} and 
\\
$\bullet$ for $p<3,$ this solution is given explicitly by the formula \fer{form:H1}. Furthermore, for $1\leq q\leq  \infty$,  if $L\in L^q(\R_+)$ then $H^{(1)}\in L^q(\R_+)$. 
\\
$\bullet$ for $p>3,$ this solution is given explicitly by the formula  \fer{form:H0}.
\end{proposition}
Because we look for an integrable function $H$ (the number of cells is supposedly finite), the function $H^{(1)}$ is preferable (take $q=1$). It also behaves better near $x \approx 0$ because the weight $p<3$ imposes that $H^{(1)}$ vanishes at $0$ as we expect.

 From the point of view of exact solutions of the direct problem, we find that $H=BN$ and $L$ belong to all spaces $ L^2(\R_+, x^p dx)$ for all $p\in \R$. Therefore, the two solutions coincide and in principle we could choose any of them. In practice, errors on the data $L$ are better handled by $H^{(1)}$ than by $H^{(2)}$ for the afore mentioned reason. 
Notice indeed that these two solutions are different in general. One can check for instance that for $L=0,$ there is a singular distributional solution $\delta_{x=0}'$. Furthermore, 
\begin{lemma}  The solutions to  \fer{eq:frag} with $L=0$ in ${\cal D}'(0,\infty)$ have the form $\f{f(\log(x))}{x^2}$ with $f \in {\cal D}' (\R)$ a $\log(2)-$ periodic distribution.
\label{lm:pbmH}
\end{lemma}

\noindent
{\bf Proof of Proposition~\ref{prop:pbmH}:} We consider the Hilbert space $X=L^2(\R_+, x^p dx)$ and we simply apply the Lax-Milgram theorem to a properly chosen bilinear form.
\\
\\ {\bf Case 1, $p<3$}. We solve the equation in the variable $y=2x$ that is (\ref{eq:BNy}).
and consider the bilinear form $a(u,v)$ on $X\times X$ defined by
$$ 
a(u,v)=4\int_0^{\infty} u(y) v(y) y^p dy -\int_0^{\infty} u(\f y 2 ) v(y) y^p dy.
$$
This form is obviously continuous and it remains to prove that it is coercive. We have 
$$
a(u,u)=4\int_0^{\infty} u(y)^2 y^p dy - \int_0^{\infty} u(\f y 2) u(y) y^p dy \geqslant (4- 2^{\f{p+1}{2}}) \; \int_0^{\infty} u (y)^2 y^p dy ,
$$
and it is indeed coercive as long as  $\alpha=4- 2^{\f{p+1}{2}} $ is positive which holds true for $p<3$.
The Lax-Milgram Theorem asserts that there is a unique $H\in X$ such that $a(H,\cdot)=(L, \cdot)$, where $( \cdot, \cdot \cdot)$ denotes the inner product in $X$, that is a solution of (\ref{eq:BNy}). 
\\
\\
{\bf Case 2, $p>3$}.  We work in the variable $x$ and consider the continuous bilinear form $b(u,v)$ on $X\times X$ defined by
$$
 b(u,v)=-4\int_0^{\infty} u(2x) v(x) x^p dx +\int_0^{\infty} u(x) v(x) x^p dx.
$$
The same calculation leads us to:
$$
b(u,u)\geqslant\alpha\int_0^{\infty} u (x)^2 x^p dx,\qquad \text{with}\quad \alpha=1- 2^{\f{3-p}{2}}>0,
$$
and the same conclusion holds. 

\

To check formulae (\ref{form:H0}) and (\ref{form:H1}), it remains to prove that these solutions belong to  the corresponding spaces:
$$
\|H^{(1)}\|_{L^2(\R_+, x^p dx)} \leqslant \sum\limits_{n=1}^\infty 2^{-2n} ||L(2^{-n} x)||_{L^2(\R_+, x^p dx)} =  \sum\limits_{n=1}^\infty 2^{\f{n}{2}(p-3)} ||L( x)||_{L^2(\R_+, x^p dx)}. 
$$
This sum converges \emph{iff} $p>3.$ In the same way, we write:
$$
||H^{(2)}||_{L^2(x^p dx)} \leqslant \sum\limits_{n=0}^\infty 2^{2n} ||L(2^n x)||_{L^2(\R_+, x^p dx)} =  \sum\limits_{n=0}^\infty 2^{\f{n}{2}(3-p)} ||L( x)||_{L^2(\R_+, x^p dx)},
$$
which converges \emph{iff} $p>3.$
\qed

\

\noindent
{\bf Proof of Lemma \ref{lm:pbmH}:} 
When $L=0,$ we first define ${\cal H}\in {\cal D}'(0,\infty)$ as the second antiderivative of $H,$ and notice that it should verify
$${\cal H}(2x)={\cal H}(x).$$
We perform the change of variables $y=\log(x)$ and notice that, if ${\cal H}\in {\cal D}'(0,\infty),$  it is equivalent to look for solutions $f \in {\cal D}'(\R)$ of 
\beq
f\big(y+\log(2)\big)=f(y).
\eeq
Hence, all the solutions in ${\cal D}' (0,\infty)$ are given by $\f{f\big(\log(x)\big)}{x^2},$ where $f\in {\cal D}'(\R).$ 
\qed

\bigskip

To conclude this Appendix, we come back to our original problem (\ref{eq:exact:inverse}) and draw the consequences in terms of $B,$ not $H.$
\begin{theorem}\label{th:exist}
Let $N \in L^2(\R_+),$ with $N (x)>0$ for $x>0.$ Let $L\in L^2(\R_+).$  
There exists a unique $B \in L^2 (\R_+,N^2 dx)$ solution of 
\beq
\label{eq:general}
4 B(2x) N(2x) - B(x) N(x)= L(x).
\eeq
\end{theorem}
\noindent
{\bf Proof:}
The theorem follows directly from Proposition~\ref{prop:pbmH} for $p=0,$ and since $N >0,$ we can  define
$B={H}/{N}$ for $B \in L^2(\R_+,N^2 dx)$.
\qed

\medskip

This theorem shows that we can find a solution $B$ of (\ref{eq:exact:inverse}) for all $N$ and all $\lb$, this is the basis of our algorithm.
However, if we want that the solution $B$ belongs to the space $L^1(\R_+;xN(x)dx)$, integration of (\ref{eq:general}) multiplied by $x$ shows that $L$ has to satisfy the condition 
$$\int_0^{\infty} x L(x)dx=0.$$
Applying this to Equation (\ref{eq:exact:inverse}), we recover that $\lb_0= {\int_0^{\infty} N(x) dx}/{\int_0^{\infty} x N(x) dx}$.
In the case of Equations (\ref{eq:invfull}) and (\ref{eq:inverse:filter}) respectively, we get formulae (\ref{eq:conserv2:PZ}) and (\ref{eq:lambdaep}), which discrete versions are expressed by (\ref{def:lbep:Q}) and (\ref{def:lbep:f}).

In view of these considerations, it is better to use a  discrete scheme defined by a matrix $A$ that preserves a similar discrete property. Namely, for all $H=(H_i),$ we should have  $\sum\limits_i i (AH)_i=0,$ in other words the vector of components $i$ belongs to the kernel of  the adjoint of $A$.
 Indeed, this property yields the (discrete) regularity $H\in L^1(\R_+;xdx).$

%

\end{document}